\documentclass[leqno,12pt]{article}

\usepackage{graphicx}
\usepackage{fancyhdr}
\usepackage{amsmath}
\usepackage{amsthm}
\usepackage{amsfonts}
\newtheorem{theorem}{Theorem}[section]
\newtheorem{lemma}[theorem]{Lemma}
\theoremstyle{definition}
\newtheorem{definition}[theorem]{Definition}

\newtheorem{notation}[theorem]{Notation}

\newcommand{\rf}[1]{(\ref{#1})}

\newcommand{\lp}[3]{|#1|_{#2}^{#3}}
\newcommand{\hp}[3]{\|#1\|_{#2}^{#3}}
\newcommand{\fr}[2]{{#1 \over #2}}

\def\N{{\Bbb N}}
\def\R{{\Bbb R}}
\def\al{{\alpha}}
\def\dl{{\delta}}
\def\O{{\Omega}}

\def\sg{{\sigma}}
\def\G{{\Gamma}}
\def\Ot{{\Omega^t}}

\def\lm{{\lambda}}

\def\iot{\intop_{\Omega^t}}

\def\vp{\varphi}

\def\vphi{{v_\varphi}}

\def\nb{\nabla}
\def\iy{\infty}

\def\divv{{\rm div}\,}
\def\curl{{\rm curl}\,}

\def\lb{\label}
\def\les{\leqslant}

\def\address#1{{\center{#1}}}
\date{}
\def\m@th{\mathsurround=0pt}
\def\eqal#1{\null\,\vcenter{\openup\jot\m@th
\ialign{\strut\hfil$\displaystyle{##}$&&$\displaystyle{{}##}$\hfil
 \crcr#1\crcr}}\,}
\def\matrix#1{\null\,\vcenter{\normalbaselines\m@th
 \ialign{\hfil$##$\hfil&&\quad\hfil$##$\hfil\crcr
 \mathstrut\crcr\noalign{\kern-\baselineskip}
 #1\crcr\mathstrut\crcr\noalign{\kern-\baselineskip}}}\,}
\def\N{{\Bbb N}}
\def\R{{\Bbb R}}
\def\divv{{\rm div}\,}
\def\curl{{\rm curl}\,}

\def\today{${\scriptscriptstyle\number\day-\number\month-\number\year}$}
\numberwithin{equation}{section}

\title{A regularity criterion for the  angular component of velocity in the norm $L_\iy(0,T;L_p(\O)),\;\fr3 p <1$ \\in axisymmetric  Navier Stokes \\equations in a cylinder }
\author{Wies\l aw J. Grygierzec$^{(1)}$, Wojciech M. Zaj\c{a}czkowski$^{(2)}$}

\begin{document}
\input amssym.def
\input amssym.tex
\maketitle
\thispagestyle{fancy}

\address{$^1$Department of Statistics and Social Policy,\\
University of Agriculture in Krak\'ow, Al. Mickiewicza 21,\\
31-120 Krak\'ow, Poland.\\
e-mail: wieslaw.grygierzec@urk.edu.pl\\
$^2$Institute of Mathematics, Polish Academy of Sciences (emeritus professor),\\
\'Sniadeckich 8, 00-656 Warsaw, Poland\\
e-mail:wz@impan.pl\\
$^2$Institute of Mathematics and Cryptology, 
Cybernetics Faculty, \\
Military University of Technology,\\
S. Kaliskiego 2, 00-908 Warsaw, Poland\\}

\begin{abstract}
We consider the axisymmetric Navier-Stokes equations in a finite cylinder $\Omega\subset\R^3$. We assume that $v_r$, $v_\varphi$, $\omega_\varphi$ vanish on the lateral part of boundary $\partial\Omega$ of the cylinder, and that $v_z$, $\omega_\varphi$, $\partial_zv_\varphi$ vanish on the top and bottom parts of the boundary $\partial\Omega$, where we used standard cylindrical coordinates, and we denoted by $\omega=\curl v$ the vorticity field. We use $H^3$ Sobolev estimates for the modified stream function (stream function divided by radius) and energy type estimates for gradient of swirl to derive two order reduction estimates. Using the estimate 
\[
\hp{\vphi}{L_\iy(0,T;L_p(\O))}{}\les A,
\]
where A is a given number and  $p>3$ we prove the existence of global regular axially-symmetric solutions.
\end{abstract}

\noindent
2020 MSC: 35A01, 35B01, 35B65, 35Q30, 76D03, 76D05\\
Key words: Navier-Stokes equations, axially-symmetric solutions, cylindrical domain

\section{Introduction}\label{s1}

We are concerned with the 3D incompressible Navier-Stokes equations,
\begin{equation}\eqal{
&\partial_tv-\nu\Delta v+v\cdot\nabla v+\nabla p=f,\cr
&\divv v=0\quad {\rm in}\ \ \Omega^T,\cr}
\label{1.1}
\end{equation}
under the axisymmetry constraint, where $\Omega^T:=\Omega\times(0,T)$, $T>0$, $v=v(x,t)\in\R^3$ denotes the velocity field, $p=p(x,t)\in\R$ denotes the pressure function, $f=f(x,t)\in\R^3$ denotes the external force field, $\nu>0$ denotes the viscosity, and $x=(x_1,x_2,x_3)$ denotes the Cartesian coordinates. As for $\Omega$ we focus on the case of a finite cylinder,
$$
\Omega=\{x\in\R^3\colon x_1^2+x_2^2<R^2,|x_3|<a\},
$$
where $a$, $R>0$ are constants. We note that
$$
S:=\partial\Omega=S_1\cup S_2,
$$
where
$$\eqal{
&S_1=\{x\in\R^3\colon\sqrt{x_1^2+x_2^2}=R,\ x_3\in[-a,a]\},\cr
&S_2=\{x\in\R^3\colon\sqrt{x_1^2+x_2^2}<R,\ x_3\in\{-a,a\}\}\cr}
$$
denote the lateral boundary and the top and bottom parts of the boundary, respectively.

In order to state the boundary conditions stating our main result we use the cylindrical coordinates $r$, $\varphi$, $z$ defined by
$$
x_1=r\cos\varphi,\quad x_2=r\sin\varphi,\quad x_3=z,
$$
and we will use standard cylindrical unit vectors, so that, for example
$$
v=v_r\bar e_r+v_\varphi\bar e_\varphi+v_z\bar e_z.
$$
We will denote partial derivatives by using the subscript comma notation, e.g.
$$
v_{r,z}:=\partial_zv_r.
$$
We assume the boundary conditions
\begin{equation}\eqal{
&v_r=v_\varphi=\omega_\varphi=0\quad &{\rm on}\ \ S_1^T=S_1\times(0,T),\cr
&v_z=\omega_\varphi=v_{\varphi,z}=0\quad &{\rm on}\ \ S_2^T=S_2\times(0,T),\cr}
\label{1.2}
\end{equation}
where $\omega:=\curl v$ denotes the vorticity vector and we assume the initial condition
\begin{equation}
v|_{t=0}=v_0,
\label{1.3}
\end{equation}
where $v_0$ is a given divergence-free vector field satisfying the same boundary conditions.

We note that such boundary conditions have first appeared in the work of Ladyzhenskaya \cite{L}. In a sense, the boundary conditions (\ref{1.2}) are natural, since, when considering the vorticity-stream function formulation we need $\omega_\varphi|_S$. This together with the no-penetration condition naturally lead to (\ref{1.2}).

We will denote the swirl by
\begin{equation}
u:=rv_\varphi.
\label{1.4}
\end{equation}
Note that
\begin{equation}\eqal{
&\omega_r=-v_{\varphi,z}=-{1\over r}u_{,z},\cr
&\omega_\varphi=v_{r,z}-v_{z,r},\cr
&\omega_z={1\over r}(rv_\varphi)_{,r}=v_{\varphi,r}+{v_\varphi\over r}={1\over r}u_{,r},\cr}
\label{1.5}
\end{equation}
so that the boundary conditions (\ref{1.2}) imply in particular that
\begin{equation}\eqal{
&\omega_r=v_{z,r}=u=0,\ \ \omega_z=v_{\varphi,r}\quad &{\rm on}\ \ S_1^T,\cr
&\omega_r=v_{r,z}=\omega_{z,z}=u_{,z}=0\quad &{\rm on}\ \ S_2^T.\cr}
\label{1.6}
\end{equation}
The Navier-Stokes equations (\ref{1.1}) in cylindrical coordinates become
\begin{equation}\eqal{
&v_{r,t}+v\cdot\nabla v_r-{v_\varphi^2\over r}-\nu\Delta v_r+\nu{v_r\over r^2}=-p_{,r}+f_r,\cr
&v_{\varphi,t}+v\cdot\nabla v_\varphi+{v_r\over r}v_\varphi-\nu\Delta v_\varphi+\nu{v_\varphi\over r^2}=f_\varphi,\cr
&v_{z,t}+v\cdot\nabla v_z-\nu\Delta v_z=-p_{,z}+f_z,\cr
&(rv_r)_{,r}+(rv_z)_{,z}=0,\cr}
\label{1.7}
\end{equation}
where
$$
v\cdot\nabla=(v_r\bar e_r+v_z\bar e_z)\cdot\nabla=v_r\partial_r+v_z\partial_z,\quad \Delta u={1\over r}(ru_{,r})_{,r}+u_{,zz}.
$$
On the other hand, the vorticity formulation becomes
\begin{equation}\eqal{
&\omega_{r,t}+v\cdot\nabla\omega_r-\nu\Delta\omega_r+\nu{\omega_r\over r^2}=\omega_rv_{r,r}+\omega_zv_{r,z}+F_r,\cr
&\omega_{\varphi,t}+v\cdot\nabla\omega_\varphi-{v_r\over r}\omega_\varphi-\nu\Delta\omega_\varphi+\nu{\omega_\varphi\over r^2}={2\over r}v_\varphi v_{\varphi,z}+F_\varphi,\cr
&\omega_{z,t}+v\cdot\nabla\omega_z-\nu\Delta\omega_z=\omega_rv_{z,r}+\omega_zv_{z,z}+F_z,\cr}
\label{1.8}
\end{equation}
where $F:=\curl f$ and the swirl is a solution to the problem
\begin{equation}\eqal{
&u_{,t}+v\cdot\nabla u-\nu\Delta u+{2\nu\over r}u_{,r}=rf_\varphi:=f_0,\cr
&u=0\quad &{\rm on}\ \ S_1^T,\cr
&u_{,z}=0\quad &{\rm on}\ \ S_2^T,\cr
&u|_{t=0}=u_0=rv_{\varphi}(0)\quad &{\rm in}\ \ \Omega.\cr}
\label{1.9}
\end{equation}
We will use the notation
\begin{equation}
(\Phi,\Gamma)=\bigg({\omega_r\over r},{\omega_\varphi\over r}\bigg),
\label{1.10}
\end{equation}
and we note that $\Phi$, $\Gamma$ satisfy
\begin{equation}
\Phi_{,t}+v\cdot\nabla\Phi-\nu\bigg(\Delta+{2\over r}\partial_r\bigg)\Phi-(\omega_r\partial_r+\omega_z\partial_z){v_r\over r}=F_r/r\equiv\bar F_r,
\label{1.11}
\end{equation}
\begin{equation}
\Gamma_{,t}+v\cdot\nabla\Gamma-\nu\bigg(\Delta+{2\over r}\partial_r\bigg)\Gamma+2{v_\varphi\over r}\Phi=F_\varphi/r\equiv\bar F_\varphi,
\label{1.12}
\end{equation}
recall [\cite{CFZ}, (\ref{1.6})]. Moreover, by (\ref{1.2}), (\ref{1.6}), $\Gamma$ and $\Phi$ satisfy the boundary conditions
\begin{equation}
\Phi=\Gamma=0\quad {\rm on}\ \ S^T.
\label{1.13}
\end{equation}
Finally, the following initial conditions are assumed
\begin{equation}
\Phi|_{t=0}=\Phi_0,\quad \Gamma|_{t=0}=\Gamma_0.
\label{1.14}
\end{equation}
We note that $(\ref{1.7})_4$ implies existence of the stream function $\psi$ which solves the problem
\begin{equation}\eqal{
&-\Delta\psi+{\psi\over r^2}=\omega_\varphi,\cr
&\psi|_S=0.\cr}
\label{1.15}
\end{equation}
Then $v$ can be expressed in terms of the stream function,
\begin{equation}\eqal{
&v_r=-\psi_{,z},\quad &v_z={1\over r}(r\psi)_{,r}=\psi_{,r}+{\psi\over r},\cr
&v_{r,r}=-\psi_{,zr},\quad &v_{z,z}=\psi_{,rz}+{\psi_{,z}\over r},\cr
&v_{r,z}=-\psi_{,zz},\quad &v_{z,r}=\psi_{,rr}+{1\over r}\psi_{,r}-{\psi\over r^2}.\cr}
\label{1.16}
\end{equation}
We will also use the modified stream function,
\begin{equation}
\psi_1:={\psi\over r},
\label{1.17}
\end{equation}

which satisfies
\begin{equation}\eqal{
&-\Delta\psi_1-{2\over r}\psi_{1,r}=\Gamma,\cr
&\psi_1|_S=0.\cr}
\label{1.18}
\end{equation}
Using the modified stream function we can express coordinates of $v$ in the form
\begin{equation}\eqal{
&v_r=-r\psi_{1,z},\quad &v_z=(r\psi_1)_{,r}+\psi_1=r\psi_{1,r}+2\psi_1,\cr
&v_{r,r}=-\psi_{1,z}-r\psi_{1,rz},\quad &v_{z,r}=3\psi_{1,r}+r\psi_{1,rr},\cr
&v_{r,z}=-r\psi_{1,zz},\quad &v_{z,z}=r\psi_{1,rz}+2\psi_{1,z}.\cr}
\label{1.19}
\end{equation}
Projecting $(\ref{1.18})_1$ on $S_2$, using $(\ref{1.18})_2$ and that $\Gamma|_{S_2}=0$ by (\ref{1.2}), we obtain
\begin{equation}
\psi_{1,zz}=0\quad {\rm on}\ \ S_2.
\label{1.20}
\end{equation}
Since in this paper we are looking for regular solutions to problem (\ref{1.1})--(\ref{1.3}), we need the following expansions near the axis of symmetry due to Liu-Wang (see \cite{LW}),
\begin{equation}\eqal{
&v_r(r,z,t)=a_1(z,t)r+a_2(z,t)r^3+\dots,\cr
&v_\varphi(r,z,t)=b_1(z,t)r+b_2(z,t)r^3+\dots,\cr
&\psi(r,z,t)=d_1(z,t)r+d_2(z,t)r^3+\dots,\cr
&\psi_1(r,z,t)=d_1(z,t)+d_2(z,t)r^2+\dots,\cr
&\psi_{1,r}(r,z,t)=2d_2(z,t)r+\dots,.\cr}
\label{1.21}
\end{equation}
In order to formulate the main results we introduce constants which depend on the initial data and forcing
\begin{notation}\label{n1.1}
$$\eqal{
&D_1=|f|_{2,1,\Omega^t}+|v(0)|_{2,\Omega}\quad &({\rm see}\ (\ref{2.5}),\cr
&D_2=|f_0|_{\infty,1,\Omega^t}+|u(0)|_{\infty,\Omega},\;f_0=rf_\vp,\;u=rv_\vp &({\rm see}\ (\ref{2.10})),\cr
&D_*=\min\{1,D_2\}\quad &({\rm see}\ (\ref{4.5})),\cr
&D_3=\frac{1}{\sqrt{2\nu}}\left(|\bar F_r|_{6/5,2,\Omega^t}+|\bar F_\varphi|_{6/5,2,\Omega^t}\right )+|\Phi(0)|_{2,\Omega}+|\Gamma(0)|_{2,\Omega}\quad &({\rm see}\ (\ref{4.7})),\cr}
$$
where $\bar F_r=F_r/r$, $\bar F_\varphi=F_\varphi/r$,
\begin{align*}
    &D_4={1\over\sqrt{\nu}}(D_1D_2+|u_{,z}(0)|_{2,\Omega}+|f_0|_{2,\Omega})\quad &({\rm see}\ (\ref{5.3})),\cr
&D_5^2=D_1^2(1+D_2)+D_1^2D_2^2+|u_{,r}(0)|_{2,\Omega}^2+|f_0|_{2,\Omega^t}^2\quad &({\rm see}\ (\ref{5.11})),\cr
&D_6^2=(D_4+D_5)\|f_\varphi\|_{L_2(0,t;L_2(S_1))}\cr
&\quad+{1\over\nu}(|F_r|_{6/5,2,\Omega^t}^2+|F_z|_{6/5,2,\Omega^t}^2)+ |\omega_r(0)|_{2,\Omega}^2+|\omega_z(0)|_{2,\Omega}^2\quad &({\rm see}\ (\ref{6.2})),\cr
&D_7=\sqrt{2}D_2^{1/2}|f_\varphi/r|_{\infty,1,\Omega^t}^{1/2}+|v_\varphi(0)|_{\infty,\Omega}\quad &({\rm see}\ (\ref{6.18})).
\end{align*}

We emphasize that the global well-posedness of axially symmetric solutions to the Navier-Stokes equations (either in the above setting or on $R^3$) remains an important open problem. We only note a few references on this topic [\cite{CFZ}, \cite{KP}, \cite{NZ}, \cite{NZ1}, \cite{NP1}, \cite{NP2}, \cite{OP}].
\end{notation}

In \cite{Z1, Z2} the second author proved the existence of global regular axially symmetric solutions by the same method as in this paper. However, he needed the following Serrin type restrictions
\begin{equation}
\psi_1|_{r=0}=0
\label{1.22}
\end{equation}
and
\begin{equation}
{|v_\varphi|_{s,\infty,\Omega^t}\over|v_\varphi|_{\infty,\Omega^t}}\ge c_0,
\label{1.23}
\end{equation}
for any $s>0$ and $c_0$ is a positive constant.

In \cite{Z1} there are assumed periodic boundary conditions on $S_2$ and in \cite{Z2} the same boundary conditions as in this paper are considered.

In \cite{OZ}, O\.za\'nski-Zaj\c{a}czkowski proved the global well-posedness assuming only condition (\ref{1.23}).

In this paper we replace restrictions (\ref{1.22}) and (\ref{1.23}), by assuming existence of a positive constant $A$ such that 
\[
\hp{\vphi}{L_\iy(0,T;L_\sg(\O))}{}\les A,\quad \sg>3\tag{1.23'}
\]
We have to recall that many regularity criterions for the angular component of velocity were found by J.Neustupa, M.Pokorny and O.Kreml (see \cite{NP1,NP2,KP}.  The restrictions were derived for the Cauchy problem and are more restrictive than (1.23'). Lately WJG. and WMZ found another criterion for regularity in \cite{GZ}. 

One of the main properties of the axisymmetric setting of the Navier-Stokes equations is the maximum principle (\ref{2.10}) for the swirl. We note that this is critical in the sence that if $v_\lambda(x,t)=\lambda v(\lambda x_1\lambda^2t)$ then $\sup|rv_{\varphi\lambda}|=\sup|rv_\varphi|$.

The main focuses of the analysis of the axisymmetric solutions to the Navier-Stokes equations are:
\begin{itemize}
\item[(1)] the best possible use of the maximum principle
\item[(2)] energy type estimates for gradient of swirl (see Lemma \ref{l5.1})
\item[(3)] the order reduction estimates for vorticity (see Lemma \ref{l6.1}).
\end{itemize}

It was demonstrated in \cite{CFZ} that the solution $v$ is controlled by the energy norm of $\Phi$, $\Gamma$,
\begin{equation}
X(t)=\|\Phi\|_{V(\Omega^t)}+\|\Gamma\|_{V(\Omega^t)},
\label{1.24}
\end{equation}
where $\|\omega\|_{V(\Omega^t)}=|\omega|_{2,\infty,\Omega^t}+|\nabla\omega|_{2,\Omega^t}$.

Now, we formulate the main result:

\begin{theorem}\label{t1.1}
Let $v$ be a smooth solution to (\ref{1.1})--(\ref{1.3}) on time interval $(0,T)$. Let the expansions (\ref{1.21}) near the axis of symmetry hold. Let the quantities from Notation \ref{n1.1} be finite for $t\in(0,T)$. Let (1.23') hold with $\sg>3$.\\
Then there exists a function $\phi=\phi(A,D_1,\dots,D_7)$ such that
\begin{equation}
X(t)\le\phi(A,D_1,\dots,D_7),
\label{1.25}
\end{equation}
where $\phi$ denotes an increasing positive function.
\end{theorem}

\begin{proof}
Lemma \ref{l4.1} yields
\begin{equation}\eqal{
X^2(t)&\le\phi(D_2,D_*,R,\delta)(1+|v_\varphi|_{\infty,\Omega^t}^{2\delta})\cdot\cr
&\quad\cdot\bigg[\intop_{\Omega^t}{v_\varphi\over r}\Phi\Gamma dxdt'+D_3^2\bigg].\cr}
\label{1.26}
\end{equation}
To prove (\ref{1.26}) we need elliptic estimates (\ref{3.14}) and (\ref{3.15}) for the modified stream function $\psi_1$. Next, we estimate in Lemma \ref{l4.2} the integral
$$
I=\intop_{\Omega^t}{v_\varphi\over r}\Phi\Gamma dxdt'.
$$
To perform it we apply the idea from \cite{CFZ}, (see Lemma 2.3)  expressing $I$ in the form
$$
I=\intop_{\Omega^t}r^dv_\varphi r^{-{1+d\over 2}}\Phi r^{-{1+d\over 2}}\Gamma dxdt'.
$$
Applying the H\"older inequality and Lemma \ref{l2.3}  we have
\begin{equation}\lb{1.27}
    |I|\les D_2^dA^{1-d}\sup_t \lp{\G}{2,\O}{\al-1/2}\lp{\Phi}{2,\Ot}{\al-1/2}\lp{\nb \Phi}{2,\Ot}{3/2-\al}\lp{\nb \G}{2,\Ot}{3/2-\al},
\end{equation}
where $\al=\fr 1 2 +\frac{(\sg-3)(1-d)}{2\sg},\;\sg>3,\;d<1.$
Using \eqref{1.27} in \rf{1.26} and Lemmas \ref{l6.1} and \ref{6.2} we prove in Theorem \ref{t4.3} the estimate \eqref{1.25}. This ends the proof.
\end{proof}

In Section \ref{s7} we prove

\begin{theorem}\label{t1.2}
Let (\ref{1.25}) hold. Let $f\in W_2^{2,1}(\Omega^t)$, $v(0)\in H^3(\Omega)$. Then
\begin{equation}
\|v\|_{W_2^{4,2}(\Omega^t)}+\|\nabla p\|_{W_2^{2,1}(\Omega^t)}\le\phi(X,\|f\|_{W_2^{2,1}(\Omega^t)},\|v(0)\|_{H^3(\Omega)}).
\label{1.33}
\end{equation}
\end{theorem}

\section{Preliminaries}\label{s2}
\subsection{Notation 2}\label{s2.1}

We will use the following notation for Lebesque spaces
$$\eqal{
&|u|_{p,\Omega}:=\|u\|_{L_p(\Omega)},\quad |u|_{p,\Omega^t}:=\|u\|_{L_p(\Omega^t)},\cr
&|u|_{p,q,\Omega^t}:=\|u\|_{L_q(0,t;L_p(\Omega))},\cr}
$$
where $p,q\in[1,\infty]$. We use standard definition of Sobolev spaces $W_p^s(\Omega)$, and we set $H^s(\Omega)=W_2^s(\Omega),\;s\in\N\cup \{0\}$, and
$$\eqal{
&\|u\|_{s,\Omega}:=\|u\|_{H^s(\Omega)},\quad &\|u\|_{s,p,\Omega}:=\|u\|_{W_p^s(\Omega)},\cr
&\|u\|_{k,p,q,\Omega^t}:=\|u\|_{L_q(0,t;W_p^k(\Omega))},\quad &\|u\|_{k,p,\Omega^t}:=\|u\|_{k,p,p,\Omega^t},\quad k\in\N\cup \{0\}.\cr}
$$
By $\phi$ we always denote an increasing positive function.

\subsection{Inequalities}\label{s2.2}

\begin{lemma} [Hardy inequality, see Lemma 2.16 in \cite{BIN}].\label{l2.1}\\
Let $p\in[1,\infty]$, $\beta\not=1/p$ and let $F(x):=\intop_0^xf(y)dy$ for $\beta>1/p$ and $F(x):=\intop_x^\infty f(y)dy$ for $\beta<1/p$. Then
\begin{equation}
|x^{-\beta}F|_{p,\R_+}\le{1\over|\beta-1/p|}|x^{-\beta+1}f|_{p,\R_+}.
\label{2.1}
\end{equation}
\end{lemma}

\begin{lemma}[Sobolev interpolation, see Sect. 15 in \cite{BIN}]\label{l2.2} 
Let $\theta$ satisfy the equality
\begin{equation}
{n\over p}-r=(1-\theta){n\over p_1}+\theta\bigg({n\over p_2}-l\bigg),\quad {r\over l}\le\theta\le 1,
\label{2.2}
\end{equation}
where $1\le p_1\le\infty$, $1\le p_2\le\infty$, $0\le r\le l$.\\
Then the interpolation holds
\begin{equation}
\sum_{|\alpha|=r}|D^\alpha f|_{p,\Omega}\le c|f|_{p_1,\Omega}^{1-\theta}\|f\|_{W_{p_2}^l(\Omega)}^\theta,
\label{2.3}
\end{equation}
where $\Omega\subset\R^n$, $D^\alpha f=\partial_{x_1}^{\alpha_1}\dots\partial_{x_n}^{\alpha_n}f$, $|\alpha|=\alpha_1+\alpha_2+\dots+\alpha_n$.
\end{lemma}

\begin{lemma}[Hardy interpolation, see Lemma 2.4 in \cite{CFZ}]\label{l2.3}
Let $0\le s\le 2$, $q\in[2,2(3-s)]$. Then there exists a positive constant $c=c(s,q)$ such that
\begin{equation}
\bigg(\intop_\Omega{|f|^{q}\over r^s}dx\bigg)^{1/q}\le c|f|_{2,\Omega}^{{3-s\over q}-{1\over 2}}|\nabla f|_{2,\Omega}^{{3\over 2}-{3-s\over q}},
\label{2.4}
\end{equation}
where $f$ does not depend on $\varphi$.
\end{lemma}

\subsection{Basic estimates}\label{s2.3}

\begin{lemma}[see Lemma 2.2 in \cite{Z1, Z2}]\label{l2.4}
Let $f\in L_{2,1}(\Omega^t)$, $v(0)\in L_2(\Omega)$. Then solutions to (\ref{1.1})--(\ref{1.3}) satisfy
\begin{equation}\eqal{
&|v(t)|_{2,\Omega}^2+\nu|\nabla v|_{2,\Omega^t}^2+\nu\intop_\Ot\bigg({v_r^2\over r^2}+{v_\varphi^2\over r^2}\bigg)dxdt'\cr
&\le 3|f|_{2,1,\Omega^t}^2+2|v(0)|_{2,\Omega}^2\equiv cD_1^2.\cr}
\label{2.5}
\end{equation}
\end{lemma}

\begin{proof}
Multiplying $(\ref{1.7})_1$ by $v_r$, $(\ref{1.7})_2$ by $v_\varphi$, $(\ref{1.7})_3$ by $v_z$, adding the results and integrating over $\Omega$ yield
\begin{equation}\eqal{
&{1\over 2}{d\over dt}\intop_\Omega(v_r^2+v_\varphi^2+v_z^2)dx+\nu\intop_\Omega(|\nabla v_r|^2+|\nabla v_\varphi|^2+|\nabla v_z|^2)dx\cr
&\quad+\nu\intop_\Omega\bigg({v_r^2\over r^2}+{v_\varphi^2\over r^2}\bigg)dx+\intop_\Omega(p_{,r}v_r+p_{,z}v_z)dx\cr
&=\intop_\Omega(f_rv_r+f_\varphi v_\varphi+f_zv_z)dx.\cr}
\label{2.6}
\end{equation}
Since $v$ is divergence free the last term on the l.h.s. vanishes. From (\ref{2.6}) we have
\begin{equation}
{d\over dt}|v|_{2,\Omega}\le|f|_{2,\Omega},
\label{2.7}
\end{equation}
where $f^2=f_r^2+f_\varphi^2+f_z^2$. Integrating (\ref{2.7}) with respect to time yields
\begin{equation}
|v(t)|_{2,\Omega}\le|f|_{2,1,\Omega^t}+|v(0)|_{2,\Omega}.
\label{2.8}
\end{equation}
Integrating (\ref{2.6}) with respect to time and using (\ref{2.8}), we obtain
\begin{equation}\eqal{
&{1\over 2}|v(t)|_{2,\Omega}^2+\nu|\nabla v|_{2,\Omega^t}^2+\nu\intop_{\Omega^t}\bigg({v_r^2\over r^2}+{v_\varphi^2\over r^2}\bigg)dxdt'\cr
&\le|f|_{2,1,\Omega^t}(|f|_{2,1,\Omega^t}+|v(0)|_{2,\Omega})+{1\over 2}|v(0)|_{2,\Omega}^2.\cr}
\label{2.9}
\end{equation}
The above inequality implies (\ref{2.5}). This concludes the proof.
\end{proof}

As for the swirl $u=rv_\varphi$, we have the following.

\begin{lemma}[Maximum principle for the swirl]\label{l2.5}
For any regular solution $v$ to (\ref{1.1})--(\ref{1.3}) we have
\begin{equation}
|u(t)|_{\infty\Omega}\le D_2=|f_0|_{\infty,1,\Omega^t}+|u(0)|_{\infty,\Omega}.
\label{2.10}
\end{equation}
\end{lemma}

\begin{proof}
Multiplying the swirl equation (\ref{1.9}) by $u|u|^{s-2}$, $s>2$, integrating over $\Omega$ and by parts, we obtain
$$
{1\over s}{d\over dt}|u|_{s,\Omega}^s+{4\nu(s-1)\over s^2}|\nabla|u|^{s/2}|_{2,\Omega}^2+{\nu\over s}\intop_\Omega(|u|^s)_{,r}drdz=\intop_\Omega f_0u|u|^{s-2}.
$$
Noting that $u|_{r=0}=u|_{r=R}=0$ (by $(\ref{1.21})_2$ and (\ref{1.2})), we see that the last term on the left-hand side vanishes, and so
$$
{d\over dt}|u|_{x,\Omega}\le|f_0|_{s,\Omega}.
$$
Integration in time and taking $s\to\infty$ gives (\ref{2.10}).
\end{proof}

\begin{lemma}[Energy estimates for $\psi$ and $\psi_1$]\label{l2.6}
For every regular solution $v$ to (\ref{1.1})--(\ref{1.3}),
\begin{equation}
\|\psi\|_{1,\Omega}^2+|\psi_1|_{2,\Omega}^2\le cD_1^2,
\label{2.11}
\end{equation}
\begin{equation}
\|\psi_{,z}\|_{1,2,\Omega^t}^2+|\psi_{1,z}|_{2,\Omega^t}^2\le cD_1^2.
\label{2.12}
\end{equation}
\end{lemma}

\begin{proof}
Multiplying $(\ref{1.15})_1$ by $\psi$, and integrating over $\Omega$ we obtain
$$\eqal{
|\nabla\psi|_{2,\Omega}^2+|\psi_1|_{2,\Omega}^2
&=\intop_\Omega\omega_\varphi\psi dx=\intop_\Omega(v_{r,z}-v_{z,r})\psi dx\cr
&=\intop_\Omega(v_z\psi_{,r}-v_r\psi_{,z})dx+\intop_{\Omega}v_z\psi_1dx\cr
&\le(|\psi_{,r}|_{2,\Omega}^2+|\psi_{,z}|_{2,\Omega}^2+|\psi_1|_{2,\Omega}^2)/2+c(v_r|_{2,\Omega}^2+|v_z|_{2,\Omega}^2),\cr}
$$
where we integrated by parts and used the boundary condition $\psi|_S=0$ (recall (\ref{1.15})) in the third equality. For (\ref{2.12}) we differentiate $(\ref{1.15})_1$ with respect to $z$, multiply by $\psi_{,z}$ and integrate over $\Omega^t$ to obtain
$$\eqal{
&\intop_{\Omega^t}|\nabla\psi_{,z}|^2dxdt'+\intop_{\Omega^t}|\psi_{1,z}|^2dxdt'= \intop_{\Omega^t}\omega_{\varphi,z}\psi_{,z}dxdt'\cr
&=-\intop_{\Omega^t}\omega_\varphi\psi_{,zz}dxdt' \le|\psi_{,zz}|_{2,\Omega^t}^2/2+c|\omega_\varphi|_{2,\Omega^t}^2,\cr}
$$
as required, where we used boundary condition $\omega_\varphi|_S=0$ (recall (\ref{1.2})) in the second equality.
\end{proof}

\section{Estimates for the modified stream function $\psi_1$}\label{s3}

Here we introduce some estimates of $\psi_1$ in terms of $\Gamma$.

\subsection{Weighted Sobolev estimates for $\psi_1$}\label{s3.1}

\begin{lemma}[see Lemma 4.2 \cite{NZ}]\label{l3.1}
Let $\mu   \in (0,1)$.  If $\psi_1$ is a sufficiently regular solution to (\ref{1.18}), then
\begin{equation}
\intop_\Omega(\psi_{1,zzz}^2+\psi_{1,rzz}^2)r^{2\mu}dx+2\mu(1-\mu)\intop_\Omega\psi_{1,zz}^2r^{2\mu-2}dx\le c\intop_\Omega\Gamma_{,r}^2r^{2\mu}dx.
\label{3.1}
\end{equation}
\end{lemma}

\begin{proof}
We differentiate (\ref{1.18}) with respect to $z$, multiply by $-\psi_{1,zzz}r^{2\mu}$ and integrate over $\Omega$ to obtain
\begin{equation}\eqal{
&\intop_\Omega\psi_{1,rrz}\psi_{1,zzz}r^{2\mu}dx+\intop_\Omega\psi_{1,zzz}^2r^{2\mu}dx\cr
&\quad+3\intop_\Omega{1\over r}\psi_{1,rz}\psi_{1,zzz}r^{2\mu}dx=-\intop_\Omega\Gamma_{,z}\psi_{1,zzz}r^{2\mu}dx.\cr}
\label{3.2}
\end{equation}
In view of (\ref{1.18}), \eqref{1.20}, the first integral on the left-hand side of (\ref{3.2}) equals
\begin{equation}\eqal{
&-\intop_\Omega\psi_{1,rrzz}\psi_{1,zz}r^{2\mu}dx=-\intop_\Omega(\psi_{1,rzz}\psi_{1,zz}r^{2\mu+1})_{,r}drdz\cr
&\quad+\intop_\Omega\psi_{1,rzz}^2r^{2\mu}dx+(2\mu+1)\intop_\Omega\psi_{1,rzz}\psi_{1,zz}r^{2\mu}drdz.\cr}
\label{3.3}
\end{equation}
Since $\psi_1|_{r=R}=\psi_{1,r}|_{r=0}=0$ (by (\ref{1.18}) and (\ref{1.21})), the first term on the right-hand side of (\ref{3.3}) vanishes. Integrating by parts with respect to $z$ in the last term on the left-hand side of (\ref{3.2}) and using (\ref{1.20}), it takes the form
\begin{equation}
-3\intop_\Omega\psi_{1,rzz}\psi_{1,zz}drdz.
\label{3.4}
\end{equation}
Using (\ref{3.3}) and (\ref{3.4}) in (\ref{3.2}) yields
\begin{equation}\eqal{
&\intop_\Omega(\psi_{1,zzz}^2+\psi_{1,rzz}^2)r^{2\mu}dx+2(\mu-1)\intop_\Omega\psi_{1,rzz}\psi_{1,zz}r^{2\mu}drdz\cr
&=-\intop_\Omega\Gamma_{,z}\psi_{1,zzz}r^{2\mu}dx.\cr}
\label{3.5}
\end{equation}
The second term on the left-hand side of (\ref{3.5}) equals
\begin{equation}\eqal{
&(\mu-1)\intop_\Omega\partial_r(\psi_{1,zz}^2)r^{2\mu}drdz\cr
&=(\mu-1)\intop_\Omega\partial_r(\psi_{1,zz}^2r^{2\mu})drdz+2\mu(1-\mu)\intop_\Omega\psi_{1,zz}^2r^{2\mu-1}drdz,\cr}
\label{3.6}
\end{equation}
where the first integral vanishes because $\psi_{1,zz}|_{r=R}=0$ (recall (\ref{1.18})) and $\psi_{1,zz}^2r^{2\mu}|_{r=0}=0$ (recall (\ref{1.21})). Using (\ref{3.6}) in (\ref{3.5}) and applying the H\"older and Young inequalities to the r.h.s. of (\ref{3.5}), we obtain (\ref{3.1}), as required.
\end{proof}

\subsection{Elliptic estimates for the modified stream function $\psi_1$}

We recall that the modified stream function $\psi_1$ is a solution to the problem (\ref{1.18}),
$$
-\Delta\psi_1-{2\over r}\psi_{1,r}=\Gamma,\quad \psi_1|_S=0.
$$
In this section we prove $H^2$ and $H^3$ elliptic estimates for $\psi_1$, in cylindrical coordinates.

\begin{lemma}[$H^2$ elliptic estimate on $\psi_1$, see Lemma 3.1 in \cite{Z1}]\label{l3.2}
If $\psi_1$ is a sufficiently regular solution to (\ref{1.18}) then
\begin{equation}\eqal{
&\intop_\Omega\bigg(\psi_{1,rr}^2+\psi_{1,rz}^2+\psi_{1,zz}^2+{\psi_{1,r}^2\over r^2}\bigg)dx\cr
&\quad+\intop_{-a}^a(\psi_{1,z}^2|_{r=0}+\psi_{1,r}^2|_{r=R})dx\le c|\Gamma|_{2,\Omega}^2.\cr}
\label{3.7}
\end{equation}
\end{lemma}

\begin{proof}
We multiply (\ref{1.18}) by $\psi_{1,zz}$ and integate over $\Omega$ to obtain
\begin{equation}
-\intop_\Omega\bigg(\psi_{1,rr}\psi_{1,zz}+\psi_{1,zz}^2+3{\psi_{1,r}\over r}\psi_{1,zz}\bigg)dx= \intop_\Omega\Gamma\psi_{1,zz}dx.
\label{3.8}
\end{equation}
Integrating by parts with respect to $r$ in the first term gives
$$\eqal{
&-\intop_\Omega(\psi_{1,r}\psi_{1,zz}r)_{,r}drdz+\intop_\Omega\psi_{1,r}\psi_{1,zzr}dx+ \intop_\Omega\psi_{1,r}\psi_{1,zz}drdz\cr
&\quad-\intop_\Omega\psi_{1,zz}^2dx-3\intop_\Omega\psi_{1,r}\psi_{1,zz}drdz= \intop_\Omega\Gamma\psi_{1,zz}dx.\cr}
$$
Thus
\begin{equation}\eqal{
&-\intop_{-a}^a[\psi_{1,r}\psi_{1,zz}r|_{r=0}^{r=R}dx+\intop_\Omega\psi_{1,r}\psi_{1,zzr}dx- \intop_\Omega\psi_{1,zz}^2dx\cr
&\quad-2\intop_\Omega\psi_{1,r}\psi_{1,zz}drdz=\intop_\Omega\Gamma\psi_{1,zz}dx.\cr}
\label{3.9}
\end{equation}
We note that the first integral vanishes since $\psi_{1,r}|_{r=0}=0$ (recall expansion (\ref{1.21})) and $\psi_{1,zz}|_{r=R}=0$. We now integrate by parts with respect to $z$ in the second and the last terms on the left-hand side and use that $\psi_{1,r}|_{S_2}=0$ (since $\psi_1|_S=0$, recall (\ref{1.18})), and we multiply by $-1$, to obtain
\begin{equation}
\intop_\Omega(\psi_{1,zr}^2+\psi_{1,zz}^2)dx-2\intop_\Omega\psi_{1,rz}\psi_{1,z}drdz=-\intop_\Omega\Gamma\psi_{1,zz}dx.
\label{3.10}
\end{equation}
We note that the last term on the left-hand side equals
$$
-\intop_\Omega(\psi_{1,z}^2)_{,r}drdz=-\intop_{-a}^a[\psi_{1,z}^2]_{r=0}^{r=R}dz= \intop_{-a}^a\psi_{1,z}^2|_{r=0}dz,
$$
since $\psi_{1,z}|_{r=R}=0$. Applying this in (\ref{3.10}), and using the Young inequality to absorb $\psi_{1,zz}$ by the left-hand side, we obtain
\begin{equation}
\intop_\Omega(\psi_{1,rz}^2+\psi_{1,zz}^2)dx+\intop_{-a}^a\psi_{1,z}^2|_{r=0}dz\le c|\Gamma|_{2,\Omega}^2.
\label{3.11}
\end{equation}
We now multiply $(\ref{1.18})_1$ by $\psi_{1,r}/r$ and integrate over $\Omega$ to obtain
\begin{equation}
3\intop_\Omega{\psi_{1,r}^2\over r^2}dx=-\intop_\Omega\bigg(\psi_{1,rr}{\psi_{1,r}\over r}+\psi_{1,zz}{\psi_{1,r}\over r}+\Gamma{\psi_{1,r}\over r}\bigg)dx.
\label{3.12}
\end{equation}
This first term on the right-hand side equals
$$
-{1\over 2}\intop_\Omega\partial_r\psi_{1,r}^2drdz=-{1\over 2}\intop_{-a}^a\psi_{1,r}^2|_{r=R}dz,
$$
where we used that $\psi_{1,r}|_{r=0}=0$ (recall expansion (\ref{1.21})) in the last equality.

As for the other terms on the right-hand side of (\ref{3.12}) we apply the Young inequality to absorb $\psi_{1,r}/r$ by the left-hand side. We obtain
$$
\intop_\Omega{\psi_{1,r}^2\over r^2}dx+{1\over 2}\intop_{-a}^a\psi_{1,r}^2|_{r=R}dz\le c(|\psi_{1,zz}|_{2,\Omega}^2+|\Gamma|_{2,\Omega}^2).
$$
The claim (\ref{3.7}) follows from this, (\ref{3.11}), and from the equation $(\ref{1.18})_1$ for $\psi_1$, which lets us estimate $\psi_{1,rr}$ in terms of $\psi_{1,zz}$, $\psi_{1,r}/r$.
\end{proof}

\begin{lemma}[$H^3$ elliptic estimates on $\psi_1$]\label{l3.3}
If $\psi_1$ is a sufficiently regular solutions to (\ref{1.18}) then
\begin{equation}
\intop_\Omega(\psi_{1,zzr}^2+\psi_{1,zzz}^2)dx+\intop_{-a}^a\psi_{1,zz}^2|_{r=0}dz\le c|\Gamma_{,z}|_{2,\Omega}^2
\label{3.13}
\end{equation}
and
\begin{equation}\eqal{
&\intop_\Omega(\psi_{1,rrz}^2+\psi_{1,rzz}^2+\psi_{1,zzz}^2)dx+\intop_{-a}^a\psi_{1,zz}^2|_{r=0}dz +\intop_{-a}^a\psi_{1,rz}^2|_{r=R}dz\cr
&\le c|\Gamma_{,z}|_{2,\Omega}^2.\cr}
\label{3.14}
\end{equation}
as well as
\begin{equation}
\bigg|{1\over r}\psi_{1,rz}\bigg|_{2,\Omega}\le c|\Gamma_{,z}|_{2,\Omega}.
\label{3.15}
\end{equation}
\end{lemma}

\begin{proof}
First we show (\ref{3.13}). We differentiate $(\ref{1.18})_1$ with respect to $z$, multiply by $-\psi_{1,zzz}$ and integrate over $\Omega$ to obtain
\begin{equation}\eqal{
&\intop_\Omega\psi_{1,rrz}\psi_{1,zzz}dx+\intop_\Omega\psi_{1,zzz}^2dx+3\intop_\Omega{1\over r}\psi_{1,rz}\psi_{1,zzz}dx\cr
&=-\intop_\Omega\Gamma_{,z}\psi_{1,zzz}dx.\cr}
\label{3.16}
\end{equation}
Integrating by parts with respect to $z$ in the first term yields
\begin{equation}
\intop_\Omega\psi_{1,rrz}\psi_{1,zzz}dx=\intop_\Omega(\psi_{1,rrz}\psi_{1,zz})_{,z}dx- \intop_\Omega\psi_{1,rrzz}\psi_{1,zz}dx,
\label{3.17}
\end{equation}
where the first term vanishes due to (\ref{1.20}). Integrating the last integral in (\ref{3.17}) by parts with respect to $r$ gives
$$
-\intop_\Omega(\psi_{1,rzz}\psi_{1,zz}r)_{,r}drdz+\intop_\Omega\psi_{1,rzz}^2dx+ \intop_\Omega\psi_{1,rzz}\psi_{1,zz}drdz,
$$
where the first integral vanishes, since $\psi_{1,rzz}|_{r=0}=\psi_{1,zz}|_{r=R}=0$ (recall (\ref{1.21}) and (\ref{1.18})). Thus, (\ref{3.16}) becomes
\begin{equation}\eqal{
&\intop_\Omega(\psi_{1,rzz}^2+\psi_{1,zzz}^2)dx+\intop_\Omega(\psi_{1,rzz}\psi_{1,zz}+ 3\psi_{1,rz}\psi_{1,zzz})drdz\cr
&=-\intop_\Omega\Gamma_{,z}\psi_{1,zzz}dx.\cr}
\label{3.18}
\end{equation}
Integrating by parts with respect to $z$ in the last term on the left-hand side of (\ref{3.18}) and using that $\psi_{1,zz}|_{S_2}=0$ (recall (\ref{1.20})) we get
\begin{equation}
\intop_\Omega(\psi_{1,rzz}+\psi_{1,zzz}^2)dx-\intop_\Omega\partial_r\psi_{1,zz}^2drdz= -\intop_\Omega\Gamma_{,z}\psi_{1,zzz}dx.
\label{3.19}
\end{equation}
Recalling (\ref{1.18}) that $\psi_{1,zz}|_{r=R}=0$ and using the Young inequality to absorb $\psi_{1,zzz}$ we obtain
$$
\intop_\Omega(\psi_{1,rzz}^2+\psi_{1,zzz}^2)dx+\intop_{-a}^a\psi_{1,zz}^2|_{r=0}dz\le c|\Gamma_{,z}|_{2,\Omega}^2,
$$
which gives (\ref{3.13}).

As for (\ref{3.14}), we differentiate $(\ref{1.18})_1$ with respect to $z$, multiply by $\psi_{1,rrz}$ and integrate over $\Omega$ to obtain
\begin{equation}
-\intop_\Omega\bigg(\psi_{1,rrz}^2+\psi_{1,zzz}\psi_{1,rrz}+3{1\over r}\psi_{1,rz}\psi_{1,rrz}\bigg)dx=\intop_\Omega\Gamma_{,z}\psi_{1,rrz}dx.
\label{3.20}
\end{equation}
We integrate the second term on the left-hand side by parts in $z$, and recall (\ref{1.20}) that $\psi_{1,zz}|_{S_2}=0$, to get
$$\eqal{
&-\intop_\Omega\psi_{1,zzz}\psi_{1,rrz}dx=\intop_\Omega\psi_{1,zz}\psi_{1,rrzz}dx\cr
&=\intop_\Omega(\psi_{1,zz}\psi_{1,rzz}r)_{,r}drdz-\intop_\Omega\psi_{1,rzz}^2dx- \intop_\Omega\psi_{1,zz}\psi_{1,rzz}drdz.\cr}
$$
We note that the first term on the right-hand side vanishes since $\psi_{1,rzz}|_{r=0}=0$ (recall (\ref{1.21})) and $\psi_{1,zz}|_{r=R}=0$ (recall (\ref{1.18})), and so (\ref{3.20}) becomes
\begin{equation}\eqal{
&\intop_\Omega(\psi_{1,rrz}^2+\psi_{1,rzz}^2)dx+\intop_\Omega(\psi_{1,zz}\psi_{1,rzz}+ 3\psi_{1,rz}\psi_{1,rrz})drdz\cr
&=-\intop_\Omega\Gamma_{,z}\psi_{1,rrz}dx.\cr}
\label{3.21}
\end{equation}
Since the first integral from the second term from the above l.h.s. equals
$$
{1\over 2}\intop_{-a}^a[\psi_{1,zz}^2]|_{r=0}^{r=R}dz=-{1\over 2}\intop_{-a}^a\psi_{1,zz}^2|_{r=0}dz
$$
(as $\psi_{1,zz}|_{r=R}=0$, recall (\ref{1.18})), and the last integral from the l.h.s. of  (\ref{3.21}) equals
$$
{3\over 2}\intop_\Omega\partial_r\psi_{1,rz}^2drdz={3\over 2}\intop_{-a}^a[\psi_{1,rz}^2]_{r=0}^{r=R}dz={3\over 2}\intop_{-a}^a\psi_{1,rz}^2|_{r=R}dz
$$
(as $\psi_{1,rz}|_{r=0}=0$, recall (\ref{1.21})). Hence (\ref{3.21}) becomes
\begin{equation}\eqal{
&\intop_\Omega(\psi_{1,rrz}^2+\psi_{1,rzz}^2)dx+\intop_{-a}^a\bigg(-{1\over 2}\psi_{1,zz}^2|_{r=0}+{3\over 2}\psi_{1,rz}^2|_{r=R}\bigg)dz\cr
&=-\intop_\Omega\Gamma_{,z}\psi_{1,rrz}dx.\cr}
\label{3.22}
\end{equation}
We now use the Young inequality to absorb $\psi_{1,rrz}$ by the left-hand side to obtain (\ref{3.14}) using also \eqref{3.13},  which in turn implies (\ref{3.15}) by differentiating $(\ref{1.18})_1$ in $z$.
\end{proof}

\section{Energy estimates for $\Phi$ and $\Gamma$}\label{s4}

\begin{lemma}\label{l4.1}
Assume that $I=\intop_\Omega{v_\varphi\over r}\Phi\Gamma dx$ be bounded for any $t\in[0,T]$. Then for regular solution to (\ref{1.1})--(\ref{1.3}), we have
\begin{equation}\eqal{
&{d\over dt}(|\Phi|_{2,\Omega}^2+|\Gamma|_{2,\Omega}^2)+\nu(|\nabla\Phi|_{2,\Omega}^2+|\nabla\Gamma|_{2,\Omega}^2)\cr
&\le{c\over D_*^2}D_2^2\cdot\bigg(1+{|v_\varphi|_{\infty,\Omega}^{2\delta}R^{2\delta}\over\delta^2D_2^{2\delta}}\bigg)\cdot\cr
&\quad\cdot\bigg(\intop_\Omega{v_\varphi\over r}\Phi\Gamma dx+{1\over 2\nu}(|\bar F_r|_{6/5,\Omega}^2+|\bar F_\varphi|_{6/5,\Omega}^2)\bigg),\cr}
\label{4.1}
\end{equation}
where $D_*=\min\{1,D_2\}$ and $D_2$ is defined in (\ref{2.10}).
\end{lemma}

\begin{proof}
Multiply (\ref{1.11}) by $\Phi$ and integrate over $\Omega$ to obtain
\begin{equation}\eqal{
&{1\over 2}{d\over dt}|\Phi|_{2,\Omega}^2+\nu|\nabla\Phi|_{2,\Omega}^2-\nu\intop_{-a}^a|\Phi^2\bigg|_{r=0}^{r=R}dz\cr
&=\intop_\Omega(\omega_r\partial_r+\omega_z\partial_z){v_r\over r}\Phi dx+\intop_\Omega\bar F_r\Phi dx,\cr}
\label{4.2}
\end{equation}
where the last term on the l.h.s. equals $\intop_{-a}^a\Phi^2|_{r=0}dz$, due to (\ref{1.6}), so it can be dropped because it is positive. Recalling (\ref{1.5}) that $\omega_r=-v_{\varphi,z}$, $\omega_z=(rv_\varphi)_{,r}/r$ we can integrate in the first term on the r.h.s. by parts,
$$\eqal{
&\intop_\Omega(\omega_r\partial_r+\omega_z\partial_z){v_r\over r}\Phi dx=\intop_\Omega\bigg(-v_{\varphi,z}\bigg({v_r\over r}\bigg)_{,r}+{(rv_\varphi)_{,r}\over r^2}v_{r,z}\bigg)\Phi dx\cr
&=\intop_\Omega v_\varphi\bigg(\bigg({v_r\over r}\bigg)_{,rz}\Phi+\bigg({v_r\over r}\bigg)_{,r}\Phi_{,z}\bigg)dx\cr
&\quad-\intop_\Omega(rv_\varphi)_{,r}\bigg({v_r\over r}\bigg)_{,z}\Phi drdz\cr
&=\intop_\Omega v_\varphi\bigg(\bigg({v_r\over r}\bigg)_{,rz}\Phi+\bigg({v_r\over r}\bigg)_{,r}\Phi_{,z}\bigg)dx\cr
&\quad-\intop_\Omega v_\varphi\bigg(\bigg({v_r\over r}\bigg)_{,zr}\Phi+\bigg({v_r\over r}\bigg)_{,z}\Phi_{,r}\bigg)dx\cr
&=-\intop_\Omega v_\varphi(\psi_{1,rz}\Phi_{,z}-\psi_{1,zz}\Phi_{,r})dx\equiv I_1,\cr}
$$
where (\ref{1.19}) was used.

Estimating, we get
$$\eqal{
|I_1|&\le\intop_\Omega\bigg|rv_\varphi{\psi_{1,rz}\over r}\Phi_{,z}\bigg|dx+\intop_\Omega \bigg|r^{1-\delta}v_\varphi{\psi_{1,zz}\over r^{1-\delta}}\Phi_{,r}\bigg|dx\cr
&\le|rv_\varphi|_{\infty,\Omega}\bigg|{\psi_{1,rz}\over r}\bigg|_{2,\Omega}|\Phi_{,z}|_{2,\Omega}+|r^{1-\delta}v_\varphi|_{\infty,\Omega}\bigg|{\psi_{1,zz}\over r^{1-\delta}}\bigg|_{2,\Omega}|\Phi_{,r}|_{2,\Omega}\cr
&\equiv I_2,\cr}
$$
where $\delta$ is assumed to be arbitrary small.

Using the maximum principle (\ref{2.10}), (\ref{3.15}), the Hardy inequality (\ref{2.1}) and (\ref{3.14}) we obtain
$$
I_2\le\varepsilon|\nabla\Phi|_{2,\Omega}^2+c(1/\varepsilon)D_2^2|\nabla\Gamma|_{2,\Omega}^2\bigg(1+ {|v_\varphi|_{\infty,\Omega}^{2\delta}R^{2\delta}\over\delta^2D_2^{2\delta}}\bigg).
$$
Using the above estimate with sufficiently small $\varepsilon$ in (\ref{4.2}) and the H\"older-Young inequalities in the last term on the r.h.s. of (\ref{4.2}) we derive
\begin{equation}\eqal{
&{1\over 2}{d\over dt}|\Phi|_{2,\Omega}^2+\nu|\nabla\Phi|_{2,\Omega}^2\le cD_2^2|\nabla\Gamma|_{2,\Omega}^2\bigg(1+ {|v_\varphi|_{\infty,\Omega}^{2\delta}R^{2\delta}\over\delta^2D_2^{2\delta}}\bigg)+c|\bar F_r|_{6/5,\Omega}^2.\cr}
\label{4.3}
\end{equation}
Multiplying (\ref{1.12}) by $\Gamma$, integrating over $\Omega$ and using the boundary conditions we have
$$\eqal{
&{1\over 2}{d\over dt}|\Gamma|_{2,\Omega}^2+\nu|\nabla\Gamma|_{2,\Omega}^2-\nu\intop_{-a}^a \Gamma^2|_{r=0}^{r=R}dz\cr
&=-2\intop_\Omega{v_\varphi\over r}\Phi\Gamma dx+\intop_\Omega\bar F_\varphi\Gamma dx.\cr}
$$
Applying the H\"older and Young inequalities to the last term and then multiplying the resulting inequality by $cD_2^2(1+|v_\varphi|_{\infty,\Omega}^{2\delta}R^{2\delta}/(\delta^2D_2^{2\delta}))$ and adding to (\ref{4.3}) gives
\begin{equation}\eqal{
&D_2^2{d\over dt}|\Gamma|_{2,\Omega}^2+\nu D_2^2|\nabla\Gamma|_{2,\Omega}^2+{d\over dt}|\Phi|_{2,\Omega}^2+\nu|\nabla\Phi|_{2,\Omega}^2\cr
&\le cD_2^2\bigg(1+{|v_\varphi|_{\infty,\Omega}^{2\delta}R^{2\delta}\over\delta^2D_2^{2\delta}}\bigg) \bigg(\intop_\Omega{v_\varphi\over r}\Phi\Gamma dx\cr
&\quad+{1\over 2\nu}(|\bar F_r|_{6/5,\Omega}^2+|\bar F_\varphi|_{6/5,\Omega}^2)\bigg).\cr}
\label{4.4}
\end{equation}
Introducing the quantity
\begin{equation}
D_*=\min\{1,D_2\}
\label{4.5}
\end{equation}
we obtain from (\ref{4.4}) the inequality (\ref{4.1}). This ends the proof.
\end{proof}

Integrating \eqref{4.1}  with respect to time yields
\begin{equation}\lb{4.6}
    \begin{aligned}
&\|\Phi\|_{V\left(\Omega^{t}\right)}^{2}+\|\Gamma\|_{V\left(\Omega^{t}\right)}^{2}\\
 &\quad\leq \frac{c}{D_{2}^{*}} D_{2}\bigg(1+\frac{\lp{\vphi}{\iy,\Ot}{2\dl}R^{2 \delta}}{\delta^{2} D_{2}^{2 \delta}}\bigg) \cdot\left(\iot \frac{\vphi}{r} \Phi \Gamma d x d t^{\prime}\right. \\
&\quad +\frac{1}{2 \nu}\left(\left|\bar F_{r}\right|_{6 / 5,2, \Omega^{t}}^{2}+\left|\bar F_{\varphi}\right|_{6 / 5,2, \Omega^{t}}^{2}\right)+|\Phi(0)|_{2, \Omega}^{2}+|\Gamma(0)|_{2, \Omega}^{2}\bigg).
\end{aligned}
\end{equation}
Simplifying, we can write (4.6) in the form
\begin{equation}\lb{4.7}
    \begin{aligned}
&\|\Phi\|_{V\left(\Omega^{t}\right)}^{2}+\|\Gamma\|_{V\left(\Omega^{t}\right)}^{2}\leqslant\left(\phi_{1}+\phi_{2}\left|\vphi\right|_{\infty, \Ot}^{2 \delta}\right)\\
 &\quad\cdot\left(\iot \frac{\vphi}{r} \Phi \Gamma d x d t^{\prime}+D_3^2\right),
\end{aligned}
\end{equation}
where $\phi_{1}=\phi_{1}\left(D_{2}\right), \phi_{2}=\phi_{2}\left(D_{2}, R, \delta\right)$.\\
Now, we estimate 
$$I=\intop_{\Omega^{t}} \frac{\vphi}{r} \Phi \Gamma d x d t^{\prime}.$$

\begin{lemma}\label{l4.2}
    Assume that there exists a positive constant $A$ such that
\begin{equation}\lb{4.8}
    \quad \sup _{t}\left|v_{\varphi}\right|_{\sigma, \Omega} \leq A
\end{equation}
 where $\sg>3$.
 Then
\begin{equation}\lb{4.9}
   |I| \leq D_{2}^{d} A^{1-d} \left|\Gamma\right|_{2, \Ot}^{\alpha-1 / 2} \lp{\Phi}{2,\Ot}{\al-1/2}\lp{\nb\Phi}{2,\Ot}{3/2-\al}\lp{\nb\G}{2,\Ot}{3/2-\al},
\end{equation}
where $\alpha=\frac{1}{2}+\frac{(\sigma-3)(1-d)}{2 \sigma}, \sigma>3,\; d<1$.
\end{lemma}
\begin{proof}
We can write $I$ in the form
\[
I=\iot \vphi r^d
r^{-\fr{1+d} 2}\Phi r^{-\fr{1+d} 2}\G dxdt'.
\]
By the Hölder inequality,
\begin{align*}
    I& \leqslant \intop_{0}^{t}\left|v_{\varphi} r^{d}\right|_{\frac{\sigma}{1-d}, \Omega}\left|r^{-\frac{1+d}{2}} \Phi\right|_{\frac{2 \sigma}{\sigma-1+d}, \Omega}\left|r^{-\frac{1+d}{2}} \Gamma\right|_{\frac{2 \sigma}{\sigma-1+d}, \Omega} d t^{\prime}\\
    &\leqslant D_{2}^{d} \intop_{0}^{t}\left|v_{\varphi}^{1-d}\right|_{\frac{\sigma}{1-d},\O}\left|r^{-\frac{1+d}{2}} \Phi\right|_{\frac{2 \sigma}{\sigma-1+d}, \Omega}\left|r^{-\frac{1+d}{2}} \Gamma\right|_{\frac{2 \sigma}{\sigma-1+d}, \Omega} d t^{\prime} \equiv I_1.
\end{align*}
The case $d=1$ yields that
$$
I_{1} \leq D_{2}\left|r^{-1} \Phi\right|_{2, \Ot}\left|r^{-1} \Gamma\right|_{2, \Ot} \equiv I_{1}^{\prime} .
$$
Since $\left.\Phi\right|_{r=0},\left.\Gamma\right|_{r=0}$ do not vanish on the axis of symmetry we can not imply the Hardy inequalities to estimate $I_{1}^{\prime}$.

Hence, we have to consider the case $d<1$.
To apply Lemma \ref{2.3} we assume
$$
s=\frac{1+d}{2} \frac{2 \sigma}{\sigma-1+d}=\frac{\sigma(1+d)}{\sigma-1+d},\; q=\frac{2 \sigma}{\sigma-1+d} .
$$
We have to check the assumptions of Lemma \ref{l2.3},
$$
s<2,\; q\in[2,2(3-s)] .
$$
The condition $s<2$ implies $\frac{\sigma(1+d)}{\sigma-1+d}<2$ so $(\sg-2)(d-1)<0$.\\
The last condition holds for $\sigma>2$ and $d<1$.\\
Next $q \leqslant 6-2 s$ implies $(\sg-3)(1-d) \geqslant 0$.\\
This holds for $\sg \geqslant 3$.\\
Applying Lemma \ref{l2.3} we have
$$
\left|r^{-\frac{1+d}{2}} \Phi\right|_{\frac{2 \sigma}{\sigma-1+d}, \Omega} \leq c|\Phi|_{2, \Omega}^{\frac{3-s}{2 q}-\frac{1}{2}}|\nabla \Phi|_{2, \Omega}^{\frac{3}{2}-\frac{3-s}{2 q}} \equiv J .
$$
We calculate
$$
\begin{aligned}
& \frac{3-s}{q}=\frac{3-\frac{\sigma(1+d)}{\sigma-1+d}}{\frac{2 \sigma}{\sigma-1+d}}=\frac{3(\sigma-1+d)-\sigma(1+d)}{2 \sigma}=\frac{1}{2}+\frac{(\sigma-3)(1-d)}{2 \sigma} \\
& \equiv \alpha .
\end{aligned}
$$
Hence $\quad \alpha-\frac{1}{2}=\frac{(\sigma-3)(1-d)}{2 \sigma}, \;\frac{3}{2}-\alpha=1-\frac{(\sigma-3)(1-d)}{2 \sigma}$.\\
Since $\alpha-1 / 2 \neq 0, \frac{3}{2}-\alpha \neq 1$ we have that $\sigma>3$.\\
Thus
$$
J \leq c|\Phi|_{2, \Omega}^{\alpha-1 / 2}|\nabla \Phi|_{2, \Omega}^{3/2-\alpha}.
$$
Using the above estimates in $I_{1}$ yields
$$
\begin{aligned}
I_{1} & \leq D_{2}^{d} \intop_{0}^{t}\left|v_{\varphi}^{1-d}\right|_{\frac{\sigma}{1-d} , \Omega}|\Phi|_{2, \Omega}^{\alpha-1 / 2}|\nabla \Phi|_{2, \Omega}^{\frac{3}{2}-\alpha}|\Gamma|_{2, \Omega}^{\alpha-\frac{1}{2}}|\nabla \Gamma|_{2, \Omega}^{\frac{3}{2}-\alpha} d t^{\prime} \\
& \equiv I_{2} .
\end{aligned}
$$
Continuing,
$$
\begin{aligned}
I_{2} & \leqslant D_{2}^{d} \sup _{t}\left|\vphi \right|_{\sigma, \Omega}^{1-d} \intop _0^{t}|\Phi|_{2, \Omega}^{\alpha-\frac{1}{2} }|\Gamma|_{2, \Omega}^{\alpha-1 / 2} \\
& \cdot |\nabla \Phi|_{2, \Omega}^{\frac{3}{2}-\alpha}|\nabla \Gamma|_{2, \Omega}^{\frac{3}{2}-\alpha} d t^{\prime} \equiv I_{3},
\end{aligned}
$$
where $\delta_{0} \in(0,1]$.\\
Applying the Hölder inequality in the time integral in $I_{3}$ we estimate it by
\begin{align*}
    &\left(\intop_{0}^{t}|\Phi|_{2, \Omega}^{(\alpha-1 / 2) \lambda_{1}} d t^{\prime}\right)^{1 / \lambda_1}\left(\intop_0^t\lp{\G}{2,\O}{(\al-1/2)\lm_2} \right)^{1/\lm_2}\left(\intop_{0}^{t}|\nabla \Phi|_{2, \Omega}^{\left(\frac{3}{2}-\alpha\right) \lambda_{3}} d t^{\prime}\right)^{1 / \lambda_{3}}\\
&\cdot\left(\intop_{0}^{t}|\nabla \Gamma|_{2, \Omega}^{\left(\frac{3}{2}-\alpha\right) \lambda_{4}} d t\right)^{1 / \lambda _4} \equiv L,
\end{align*}
where
$$
1 / \lambda_{1}+1 / \lambda_{2}+1 / \lambda_{3}+1/\lm_4=1
$$
Assuming
$(\alpha-1 / 2)  \lambda_{1}=2,\;(\alpha-1 / 2)  \lambda_{2}=2,\;\left(\frac{3}{2}-\alpha\right) \lambda_{3}=2,\;\left(\frac{3}{2}-\alpha\right) \lambda_{4}=2$ we derive the restriction
$$
\left(\alpha-\frac{1}{2}\right) +\frac{3}{2}-\alpha=1 \quad \text { so } \quad\left(\alpha-\frac{1}{2}\right)=\alpha-\frac{1}{2}.
$$
Hence 
$$L \leq|\Phi|_{2, \Ot}^{\alpha-1 / 2}\lp{\G}{2,\Ot}{\al-1/2}|\nabla \Phi|_{2, \Ot}^{\frac{3}{2}-\alpha}|\nabla \Gamma|_{2, \Ot}^{\frac{3}{2}-\alpha}.$$
Using the estimate in $I_{3}$ yields
$$I_{3} \leq D_{2}^{d} \sup _{t}\left|v_{\varphi}\right|_{\sigma, \Omega}^{1-d} |\Gamma|_{2, \Ot}^{\alpha-1 / 2}|\Phi|_{2, \Ot}^{\alpha-1 / 2}|\nabla \Phi|_{2, \Ot}^{\frac{3}{2}-\alpha}|\nabla \Gamma|_{2, \Ot}^{\frac{3}{2}-\alpha}.$$
The above estimate implies \eqref{4.9} and concludes the proof.
\end{proof}
Introduce the notation
\begin{equation}\lb{4.10}
     X^{2}(t)=\|\Phi\|_{V\left(\Omega^{t}\right)}^{2}+\|\Gamma\|_{V\left(\Omega^{t}\right)}^{2}.
\end{equation}
Applying Lemmas \ref{l4.1}, \ref{l4.2} we have
\begin{theorem}\label{t4.3}
Assume that $\sigma>3,\; d<1$, and for a given constant $A$ we have
$$
\sup _{t}\left|v_{\varphi}\right|_{\sigma, \Omega} \leqslant A .
$$
Assume that all quantities from Notation 1.1  are finite.\\
Then
\begin{equation}\lb{4.11}
   X(t) \leq \phi\left(A, D_{1}, \cdots, D_{7}, R, \delta\right).  \end{equation}
    
\end{theorem}
\begin{proof}
    Using notation \eqref{4.9} and \eqref{4.10}  in \eqref{4.7}  yields
\begin{equation}\lb{4.12}
    X^{2} \leq\left(\phi_{1}+\phi_{2}\left|\vphi\right|_{\iy, \Ot}^{2 \delta}\right)\left(D_2^dA^{1-d}\lp{\Phi}{2,\Ot}{\al-1/2}X^{5/2-\al}\right)+D_3^2.
\end{equation}    
From \eqref{6.17} we have
\begin{equation}\lb{4.13}
    \left|v_{\varphi}\right|_{\infty, \Ot} \leqslant \frac{D_{2}}{\sqrt{\nu}} D_{1}^{1 / 4} X^{3 / 4}+D_{7}
\end{equation}
and \eqref{6.1} implies
\begin{equation}\lb{4.14}
    \begin{split}
        |\Phi|_{2, \Omega t}^{2} \leqslant \phi_{3}\left(D_{1}, D_{2}, D_{4}, D_{5}, R, \delta\right)\left( 1+\lp{\vphi}{2,\Ot}{2\dl}\right)|\nabla\G|_{2, \Ot}+c D_{6}^{2}.\\
    \end{split}
\end{equation}

Using \eqref{4.13} and \eqref{4.14} in \eqref{4.12} yields
\begin{equation}\lb{4.15}
    \begin{aligned}
        X^{2} &\leq\left[\phi_{1}+\phi_{2}\left(\frac{D_{2}}{\sqrt{\nu}} D_{1}^{1 / 4} X^{3 / 4}+D_{7}\right)^{2 \delta}\right]\\
       & \cdot\left[D_{2}^{d} A^{1-\alpha}\left(\phi_{3}^{1 / 2}\left(1+\left(\frac{D_{2}}{\sqrt{\nu}} D_{1}^{1 / 4} X^{3 / 4}+D_{7}\right)^{\delta}\right) X+D_{6}\right)^{\frac{\alpha}{2}-\frac{1}{4}} \right. \\
&\cdot X^{\frac{5}{2}-\alpha} +D_{3}^{2}\bigg].
    \end{aligned}
\end{equation}
Simplifying, we can write \eqref{4.15} in the form
\begin{equation}\lb{4.16}
    X^{2} \leq \phi_{4} X^{\frac{3}{2} \delta+\frac{3}{4}\left(\frac{\alpha}{2}-\frac{1}{4}\right) \delta} X^{\frac{5}{2}-\alpha+\frac{\alpha}{2}-\frac{1}{4}}+\phi_{5},
\end{equation}
where $\phi_{i}=\phi_{i}\left(A, D_{1}, D_{2}, D_{3}, D_{4}, D_{5}, D_{6}, R, \delta\right), i=4,5$,\\
and
$$
\begin{aligned}
& \frac{5}{2}-\alpha+\frac{\alpha}{2}-\frac{1}{4}=\frac{9}{4}-\frac{\alpha}{2}=2-\left(\frac{\alpha}{2}-\frac{1}{4}\right)=2-\frac{(\sigma-3)(1-d)}{2 \sigma} \\
& \equiv 2-\alpha_{0}.
\end{aligned}
$$
Then \eqref{4.16} takes the form
\begin{equation}\lb{4.17}
    X^{2} \leq \phi_{4} X^{2-\alpha_{0}+\left(\frac{3}{2}+\frac{3}{4} \alpha_{0}\right) \delta}+\phi_{5}.
\end{equation}
Since $\dl$ is small we have
$$
\alpha_{0}-\frac{3}{2}\left(1+\frac{\alpha_{0}}{2}\right) \delta>0 .
$$

Then by the Young inequality we derive \eqref{4.11}. This ends the proof,

\end{proof}

\section{Estimates for swirl $u=rv_\varphi$}\label{s5}

We derive energy estimate for $\nabla u$. Recall that swirl $u=rv_\varphi$ satisfies
\begin{equation}\eqal{
&u_{,t}+v\cdot\nabla u-\nu\Delta u+{2\nu\over r}u_{,r}=f_0,\cr
&u=0\quad &{\rm on}\ \ S_1,\cr
&u_{,z}=0\quad &{\rm on}\ \ S_2.\cr}
\label{5.1}
\end{equation}

\begin{lemma}[see Lemma 5.1 in \cite{Z1,Z2,OZ}]\label{l5.1}
Any regular solution $u$ to (\ref{5.1}) satisfies
\begin{equation}
|u_{,z}(t)|_{2,\Omega}^2+\nu|\nabla u_{,z}|_{2,\Omega^t}^2\le cD_4^2,
\label{5.2}
\end{equation}
\begin{equation}
|u_{,r}(t)|_{2,\Omega}^2+\nu(|u_{,rr}|_{2,\Omega^t}^2+|u_{,rz}|_{2,\Omega^t}^2)\le cD_5^2,
\label{5.3}
\end{equation}
where $D_4^2={1\over\nu}(D_1^2D_2^2+|u_{,z}(0)|_{2,\Omega}^2+|f_0|_{2,\Omega}^2)$ and $D_5$ is defined in (\ref{5.11}).
\end{lemma}

\begin{proof}
Differentiating (\ref{5.1}) with respect to $z$, multiplying by $u_{,z}$ and integrating over $\Omega$ we obtain
\begin{equation}
\begin{aligned}
    &{1\over 2}{d\over dt}|u_{,z}|_{2,\Omega}^2-\nu\intop_\Omega\divv(\nabla u_{,z}u_{,z})dx+\nu\intop_\Omega|\nabla u_{,z}|^2dx\cr
&+2\nu\intop_\Omega u_{,zr}u_{,z}drdz+\intop_\Omega v_{,z}\cdot\nabla uu_{,z}dx+{1\over 2}\intop_\Omega v\cdot\nabla(u_{,z}^2)dx\cr
&=\intop_\Omega f_{0,z}u_{,z}dx.\cr
\end{aligned}
\label{5.4}
\end{equation}
The second term vanishes due to the boundary condition $u_{,z}|_S=0$ (recall $(\ref{1.2})_{1,2}$. The fourth term equals
$$
\nu\intop_\Omega\partial_r(u_{,z}^2)drdz=\nu\intop_{-a}^au_{,z}^2|_{r=0}^{r=R}dz=0,
$$
because $u_{,z}|_{r=R}=0$ (see $(\ref{1.2})_{1,2}$) and $u_{,z}|_{r=0}=0$ (recall (\ref{1.21})). Similarly, the sixth term equals
$$
{1\over 2}\intop_\Omega v\cdot\nabla u_{,z}^2dx={1\over 2}\intop_Sv\cdot\bar nu_{,z}^2dS=0,
$$
because $v\cdot\bar n|_S=0$ (recall $(\ref{1.2})_{1,2}$). Integrating by parts in the fifth term in (\ref{5.4}), and noting that the boundary term vanishes (since $u_{,z}=0$ on $S$), we obtain
$$
\bigg|\intop_\Omega(v_{,z}\cdot\nabla)u\cdot u_{,z}dx\bigg|=\bigg|\intop_\Omega v_{,z}\cdot\nabla u_{,z}udx\bigg|\le\delta\intop_\Omega|\nabla u_{,z}|^2dx+{c\over\delta}|u|_{\infty,\Omega}^2\intop_\Omega v_{,z}^2dx.
$$
Finally, integrating the right-hand side of (\ref{5.4}) by parts in $z$ we obtain
$$
\intop_\Omega f_{0,z}u_{,z}dx=-\intop_\Omega f_0u_{,zz}dx\le\delta|u_{,zz}|_{2,\Omega}^2+{c\over\delta}|f_0|_{2,\Omega}^2,
$$
where we used that $\intop_{S_2}|f_0u_{,z}|_{z=-a}^{z=a}rdr=0$ because $u_{,z}|_{S_2}=0$ (recall $(\ref{1.2})_2$). Using the above results in (\ref{5.4}) gives
$$
{d\over dt}|u_{,z}|_{2,\Omega}^2+\nu|\nabla u_{,z}|_{2\Omega}^2\le{c\over\nu}(|u|_{\infty,\Omega}^2|v_{,z}|_{2,\Omega}^2+|f_0|_{2,\Omega}^2).
$$
Integrating in $t\in(0,T)$ gives
\begin{equation}\eqal{
&|u_{,z}(t)|_{2,\Omega}^2+\nu|\nabla u_{,z}|_{2,\Omega^t}^2\le{c\over\nu}|u|_{\infty,\Omega^t}^2|v_{,z}|_{2,\Omega^t}^2\cr
&\quad+|u_{,z}(0)|_{2,\Omega}^2+|f_0|_{2,\Omega^t}^2)\le cD_4^2\cr}
\label{5.5}
\end{equation}
which proves (\ref{5.2}) because energy estimate (\ref{2.5}) and maximum principle (\ref{2.10}) for the swirl $u$ were used. To prove (\ref{5.3}) we differentiate $(\ref{5.1})_1$ with respect to $r$, multiply the resulting equation by $u_{,r}$ and integrate over $\Omega$ to obtain
\begin{equation}\eqal{
&{1\over 2}{d\over dt}|u_{,r}|_{2,\Omega}^2+\intop_\Omega v_{,r}\cdot\nabla uu_{,r}dx+\intop_\Omega v\cdot\nabla u_{,r}u_{,r}dx\cr
&\quad-\nu\intop_\Omega(\Delta u)_ru_{,r}dx+2\nu\intop_\Omega u_{,rr}u_{,r}drdz-2\nu\intop_\Omega{u_{,r}^2\over r^2}dx\cr
&=\intop_\Omega f_{0,r}u_{,r}dx.\cr}
\label{5.6}
\end{equation}
We now examine particular terms in (\ref{5.6}). The second term equals
$$\eqal{
&\intop_\Omega v_{,r}\cdot\nabla uu_{,r}rdrdz=\intop_\Omega(rv_{r,r}u_{,r}+rv_{z,r}u_{,z})u_{,r}drdz\cr
&=-\intop_\Omega[(rv_{r,r}u_{,r})_{,r}+(rv_{z,r}u_{,r})_{,z}]udrdz\equiv I,\cr}
$$
where we integrated by parts with respect to $r$ and $z$, respectively, and used the boundary conditions $u|_{S_1}=u|_{r=0}=0$ (recall (\ref{1.2}) and (\ref{1.21})) and $v_{z,r}|_{S_2}=0$ (recall (\ref{1.2})). Continuing, we have
$$\eqal{
I&=-\intop_\Omega[(rv_{r,r})_{,r}+(rv_{z,r})_{,z}]u_{,r}udrdz\cr
&\quad-\intop_\Omega[rv_{r,r}u_{,rr}+rv_{z,r}u_{,rz}]udrdz=I_1+I_2.\cr}
$$
Differentiating the divergence-free equation $(\ref{1.7})_4$ in $r$ gives $v_{r,rr}+v_{z,zr}+{v_{r,r}\over r}-{v_r\over r^2}=0$. Hence $I_1$ equals
$$
I_1=-\intop_\Omega{v_r\over r}u_{,r}udrdz.
$$
Using the Young inequality in $I_2$ yields
$$
|I_2|\le{\nu\over 2}(|u_{,rr}|_{2,\Omega}^2+ |u_{,rz}|_{2,\Omega}^2)+{c\over\nu}|u|_{\infty,\Omega}^2(|v_{r,r}|_{2,\Omega}^2+|v_{z,r}|_{2,\Omega}^2).
$$
The third term on the l.h.s. of (\ref{5.6}) equals
$$
{1\over 2}\intop_\Omega v\cdot\nabla u_{,r}^2dx={1\over 2}\intop_\Omega\divv(vu_{,r}^2)dx=0,
$$
since $v\cdot\bar n|_S=0$. As for the fourth term in (\ref{5.6}) we have
$$\eqal{
&-\intop_\Omega(\Delta u)_{,r}u_{,r}dx=-\intop_\Omega\bigg(u_{,rrr}+\bigg({1\over r}u_{,r}\bigg)_{,r} +u_{,rzz}\bigg)u_{,r}rdrdz\cr
&=-\intop_\Omega\bigg[\bigg(u_{,rr}+{1\over r}u_{,r}\bigg)u_{,r}r\bigg]_{,r}drdz+\intop_\Omega u_{,rr}(u_{,r}r)_{,r}drdz\cr
&\quad+\intop_\Omega{1\over r}u_{,r}(u_{,r}r)_{,r}drdz+\intop_\Omega u_{,rz}^2dx\cr
&=-\intop_{-a}^a\bigg[\bigg(u_{,rr}+{1\over r}u_{,r}\bigg)u_{,r}r\bigg]\bigg|_{r=0}^{r=R}dz+ \intop_\Omega(u_{,rr}^2+u_{,rz}^2)dx\cr
&\quad+\intop_\Omega{u_{,r}^2\over r^2}dx+2\intop_\Omega u_{,rr}u_{,r}drdz.\cr}
$$
Using the above expressions in (\ref{5.6}) yields
\begin{equation}\eqal{
&{1\over 2}{d\over dt}|u_{,r}|_{2,\Omega}^2+{\nu\over 2}\intop_\Omega(u_{,rr}^2+u_{,rz}^2)dx-\nu\intop_\Omega{u_{,r}^2\over r^2}dx\cr
&\quad-\nu\intop_{-a}^a\bigg[\bigg(u_{,rr}+{1\over r}u_{,r}\bigg)u_{,r}r\bigg]\bigg|_{r=0}^{R=R}dz+4\nu\intop_\Omega u_{,rr}u_{,r}drdz\cr
&\le\intop_\Omega f_{0,r}u_{,r}dx+\intop_\Omega{v_r\over r}{u_{,r}\over r}udx+{c\over\nu}|u|_{\infty,\Omega}^2(|v_{r,r}|_{2,\Omega}^2+|v_{z,r}|_{2,\Omega}^2).\cr}
\label{5.7}
\end{equation}
The last term on the l.h.s. of (\ref{5.7}) equals
$$
2\nu\intop_{-a}^au_{,r}^2\bigg|_{r=0}^{r=R}dz=2\nu\intop_{-a}^au_{,r}^2\bigg|_{r=R}dz.
$$
Since the expansion (\ref{1.21}) implies that
\begin{equation}
u=b_1(z,t)r^2+b_2(z,t)r^3+\dots,
\label{5.8}
\end{equation}
thus that $u_{,r}|_{r=0}=0$.

Projecting $(\ref{5.1})_1$ on $S_1$ implies
\begin{equation}\eqal{
&-\nu\intop_{-a}^a\bigg(u_{,rr}+{1\over r}u_{,r}\bigg)u_{,r}r\bigg|_{r=R}dz=-2\nu\intop_{-a}^au_{,r}^2\bigg|_{r=R}dz\cr
&\quad+\intop_{-a}^af_0u_{,r}r\bigg|_{r=R}dz,\cr}
\label{5.9}
\end{equation}
where in the last equality we used $(\ref{1.9})_1$ projected onto $S_1$. Integration by parts in $r$ in the first term on the r.h.s. of (\ref{5.7}) gives
$$
\intop_{-a}f_0u_{,r}r|_{r=R}dz-\intop_\Omega f_0u_{,rr}dx-\intop_\Omega f_0u_{,r}drdz,
$$
where we used (\ref{5.8}) again to note that $u_{,r}|_{r=0}=0$. We note again that the first term above cancells with the last term of (\ref{5.9}), while the remaining terms can be estimated using the Young inequality by
$$
{\nu\over 4}|u_{,rr}|_{2,\Omega}^2+\nu\bigg|{u_{,r}\over r}\bigg|_{2,\Omega}^2+{c\over\nu}|f_0|_{2,\Omega}^2.
$$
Using the above estimates in (\ref{5.7}) and simplifying we get
\begin{equation}\eqal{
&{1\over 2}{d\over dt}|u_{,r}|_{2,\Omega}^2+{\nu\over 4}(|u_{,rr}|_{2,\Omega}^2+|u_{,rz}|_{2,\Omega}^2)\le 2\nu\intop_\Omega{u_{,r}^2\over r^2}dx\cr
&\quad+\intop_\Omega{v_r\over r}{u_{,r}\over r}udx+c|u|_{\infty,\Omega}^2(|v_{r,r}|_{2,\Omega}^2+|v_{z,r}|_{2,\Omega}^2)+{c\over\nu}|f_0|_{2,\Omega}^2.\cr}
\label{5.10}
\end{equation}
Integrating (\ref{5.10}) with respect to time yields
\begin{equation}\eqal{
&|u_{,r}(t)|_{2,\Omega}^2+\nu(|u_{,rr}|_{2,\Omega^t}^2+|u_{,rz}|_{2,\Omega^t}^2)\cr
&\le cD_1^2(1+D_2)+cD_1^2D_2^2+|f_0|_{2,\Omega^t}^2\cr
&\quad+|u_{,r}(0)|_{2,\Omega}^2\equiv cD_5^2.\cr}
\label{5.11}
\end{equation}
This implies (\ref{5.3}) and concludes the proof.
\end{proof}

\section{Order reduction estimates}\label{s6}

\begin{lemma}\label{l6.1}
(See \cite{Z1,Z2,OZ}) Any regular solution to (\ref{1.1})--(\ref{1.3}) satisfies
\begin{equation}\eqal{
&\|\omega_r\|_{V(\Omega^t)}^2+\|\omega_z\|_{V(\Omega^t)}^2+\bigg|{\omega_r\over r}\bigg|_{2,\Omega^t}^2\cr
&\le{1\over\nu}\phi(D_1,D_2,D_4,D_5)\bigg({R^{\varepsilon_0}\over\varepsilon_0} |v_\phi|_{\infty,\Omega^t}^{\varepsilon_0}+{R^{2\varepsilon_0}\over\varepsilon_0^2} |v_\phi|_{\infty,\Omega^t}^{2\varepsilon_0}\bigg)|\nabla\Gamma|_{2,\Omega^t}\cr
&\quad+cD_6^2,\cr}
\label{6.1}
\end{equation}
where
\begin{equation}\eqal{
D_6^2&=(D_4+D_5)\|f_\varphi\|_{L_2(0,t;L_3(S_1)}\cr
&\quad+{1\over\nu}(|F_r|_{6/5,2/\Omega^t}^2+|F_z|_{6/5,2/\Omega^t}^2)+|\omega_r(0)|_{2,\Omega}^2+ |\omega_z(0)|_{2,\Omega}^2.\cr}
\label{6.2}
\end{equation}
\end{lemma}

\begin{proof}
Multiplying $(\ref{1.8})_1$ by $\omega_r$, $(\ref{1.8})_3$ by $\omega_z$, adding the resulting equations and integrating over $\Omega^t$, we obtain
\begin{equation}\eqal{
&{1\over 2}(|\omega_r(t)|_{2,\Omega}^2+|\omega_z(t)|_{2,\Omega}^2)\cr
&\quad+\nu\bigg(|\nabla\omega_r|_{2,\Omega^t}^2+|\nabla\omega_z|_{2,\Omega^t}^2+\bigg|{\omega_r\over r}\bigg|_{2,\Omega^t}^2\bigg)\cr
&=\nu\intop_{S^t}(\bar n\cdot\nabla\omega_z\omega_z+\bar n\cdot\nabla\omega_r\omega_r)dSdt'\cr
&\quad+\intop_{\Omega^t}(v_{r,r}\omega_r^2+v_{z,z}\omega_z^2+(v_{r,z}+v_{z,r})\omega_r\omega_z)dxdt'\cr
&\quad+\intop_{\Omega^t}(F_r\omega_r+F_z\omega_zdxdt'\cr
&\quad+{1\over 2}(|\omega_r(0)|_{2,\Omega}^2+|\omega_z(0)|_{2,\Omega}^2)\equiv I_1+J+I_2\cr
&\quad+{1\over 2}(|\omega_r(0)|_{2,\Omega}^2+|\omega_z(0)|_{2,\Omega}^2).\cr}
\label{6.3}
\end{equation}
First we examine $I_1$. Since $\omega_r=-v_{\varphi,z},v_\varphi|_{r=R}=0$ and $v_{\varphi,z}|_{S_2}=0$ we obtain 
$$
\intop_S\bar n\cdot\nabla\omega_r\omega_rdS=0.
$$
Using $(\ref{1.5})_3$ we get $\omega_z=v_{\varphi,r}+{v_\varphi\over r}$. Since $v_{\varphi,z}|_{S_2}=0$ we have
$$\eqal{
&-\nu\intop_{S^t}\bar n\cdot\nabla\omega_z\omega_zdS_1dt'=-\nu\intop_{S_1^t}\bar n\cdot\nabla\omega_z\omega_zdS_1dt'\cr
&=-\nu\intop_{S_1^t}\partial_r\bigg(v_{\varphi,r}+{v_\varphi\over r}\bigg)\bigg(v_{\varphi,r}+{v_\varphi\over r}\bigg)Rdzdt'\cr
&=-\nu\intop_0^t\intop_{-a}^a\bigg(v_{\varphi,rr}+{v_{\varphi,r}\over r}\bigg)v_{\varphi,r}\bigg|_{r=R}Rdzdt'\equiv I_1^1,\cr}
$$
where we used that $v_\varphi|_{S_1}=0$. Projecting $(\ref{1.7})_2$ on $S_1$ yields
$$
-\nu\bigg(v_{\varphi,rr}+{1\over r}v_{\varphi,r}\bigg)=f_\varphi.
$$
Hence
$$
I_1^1=R\intop_0^t\intop_{-a}^af_\varphi v_{\varphi,r}|_{r=R}dzdt'=\intop_0^t\intop_{-a}^af_\varphi\bigg(u_{,r}-{1\over R}u\bigg)dzdt'
$$
and
$$
|I_1^1|\le|f_\varphi|_{2,S_1^t}(|u_{,r}|_{2,S_1^t}+|u|_{2,S_1^t})\le c|f_\varphi|_{2,S_1^t}(D_4+D_5).
$$
Summarizing,
\begin{equation}
I_1\le c|f_\varphi|_{2,S_1^t}(D_4+D_5).
\label{6.4}
\end{equation}
Next, we examine $I_2$. By the H\"older inequality, we get
\begin{equation}\eqal{
I_2&\le\varepsilon(|\omega_r|_{6,2,\Omega^t}^2+|\omega_z|_{6,2,\Omega^t}^2)\cr
&\quad+{1\over 4\varepsilon}(|F_r|_{6/5,2,\Omega^t}^2+|F_z|_{6/5,2,\Omega^t}^2).\cr}
\label{6.5}
\end{equation}
Finally, we examine
\begin{equation}
J=\intop_{\Omega^t}[v_{r,r}\omega_r^2+v_{z,z}\omega_z^2+(v_{r,z}+v_{z,r})\omega_r\omega_z]dxdt'.
\label{6.6}
\end{equation}
Using (\ref{1.5}) and (\ref{1.16}) yields
\begin{equation}\eqal{
J&=\intop_{\Omega^t}\bigg[-\psi_{,zr}\bigg({1\over r}u_{,z}\bigg)^2+\bigg(\psi_{,rz}+{\psi_{,z}\over r}\bigg)\bigg({1\over r}u_{,r}\bigg)^2\cr
&\quad-\bigg(-\psi_{,zz}+\psi_{,rr}+{1\over r}\psi_{,r}-{\psi\over r^2}\bigg)\bigg({1\over r}u_{,z}\bigg)\bigg({1\over r}u_{,r}\bigg)\bigg]dxdt'\cr}
\label{6.7}
\end{equation}
Consider $J_1$. Integrating by parts with respect to $z$ and using that $u_{,z}|_{S_2}=0$, we obtain
$$\eqal{
J_1&=-\intop_{\Omega^t}\psi_{,zr}{1\over r}u_{,z}{1\over r}u_{,z}dxdt'=\intop_{\Omega^t}\psi_{,zzr}{1\over r^2}u_{,z}udxdt'\cr
&\quad+\intop_{\Omega^t}\psi_{,zr}{1\over r^2}u_{,zz}udxdt'\equiv J_{11}+J_{12}.\cr}
$$
Using the transformation $\psi=r\psi_1$, we have
$$
J_{11}=\intop_{\Omega^t}\bigg({\psi_{1,zz}\over r}+\psi_{1,zzr}\bigg){u_{,z}\over r}udxdt'\equiv J_{11}^1+J_{11}^2.
$$
By the H\"older inequality,
$$\eqal{
|J_{11}^1|&\le\intop_{\Omega^t}\bigg|{\psi_{1,zz}\over r^{1-\varepsilon_0}}\bigg|\,\bigg|{u_{,z}\over r}\bigg|\,|v_\varphi|^{\varepsilon_0}|u|^{1-\varepsilon_0}dxdt'\cr
&\le D_2^{1-\varepsilon_0}|v_\varphi|_{\infty,\Omega^t}^{\varepsilon_0}\bigg|{u_{,z}\over r}\bigg|_{2,\Omega^t}\bigg|{\psi_{1,zz}\over r^{1-\varepsilon_0}}\bigg|_{2,\Omega^t},\cr}
$$
where (\ref{2.10}) is used. In view of (\ref{2.5})
$$
\bigg|{u_{,z}\over r}\bigg|_{2,\Omega^t}\le|v_{\varphi,z}|_{2,\Omega^t}\le D_1
$$
and (\ref{2.1}), (\ref{3.13}) imply
$$
\bigg|{\psi_{1,zz}\over r^{1-\varepsilon_0}}\bigg|\le c{R^{\varepsilon_0}\over\varepsilon_0}|\psi_{1,zzr}|_{2,\Omega^t}\le c{R^{\varepsilon_0}\over\varepsilon_0}|\Gamma_{,z}|_{2,\Omega^t}.
$$
Summarizing,
$$
|J_{11}^1|\le c{R^{\varepsilon_0}\over\varepsilon_0}D_1D_2^{1-\varepsilon_0}|v_\varphi|_{\infty,\Omega^t}^{\varepsilon_0} |\Gamma_{,z}|_{2,\Omega^t}.
$$
Next,
$$
|J_{11}^2|\le|u|_{\infty,\Omega^2}\bigg|{u_{,z}\over r}\bigg|_{2,\Omega^t}|\psi_{1,zzr}|_{2,\Omega^t}\le cD_2D_1|\Gamma_{,z}|_{2,\Omega^t},
$$
where (\ref{2.5}), (\ref{2.10}), (\ref{3.13}) were used. Hence,
\begin{equation}
|J_{11}|\le cD_1D_2^{1-\varepsilon_0}{R^{\varepsilon_0}\over\varepsilon_0} |v_\varphi|_{\infty,\Omega^t}^{\varepsilon_0}|\Gamma_{,z}|_{2,\Omega^t}+cD_1D_2|\Gamma_{,z}|_{2,\Omega^t}.
\label{6.8}
\end{equation}
Next,
$$
J_{12}=\intop_{\Omega^t}\bigg({\psi_{1,z}\over r^2}+{\psi_{1,zr}\over r}\bigg)u_{,zz}udxdt\equiv J_{12}^1+J_{12}^2.
$$
Estimates (\ref{2.10}), (\ref{3.15}), (\ref{5.1}) imply
$$
|J_{12}^2|\le|u|_{\infty,\Omega^t}|u_{,zz}|_{2,\Omega^t}\bigg|{\psi_{1,zr}\over r}\bigg|_{2,\Omega^t}\le cD_2D_4|\Gamma_{,z}|_{2,\Omega^t}.
$$
Hence,
\begin{equation}
|J_{12}|\le\bigg|\intop_{\Omega^t}{\psi_{1,z}\over r^2}u_{,zz}udxdt'\bigg|+cD_2D_4|\Gamma_{,z}|_{2,\Omega^t}.
\label{6.9}
\end{equation}
Definition of $J_1$ and (\ref{6.8}), (\ref{6.9}) imply
\begin{equation}\eqal{
|J_1|&\le c\bigg[{R^{\varepsilon_0}\over\varepsilon_o}D_1D_2^{1-\varepsilon_0}|v_\varphi|_{\infty,\Omega^t}^{\varepsilon_0} +D_1D_2+D_2D_4\bigg]|\Gamma_{,z}|_{2,\Omega^t}\cr
&\quad+\bigg|\intop_{\Omega^t}{\psi_{1,z}\over r^2}u_{,zz}udxdt\bigg|.\cr}
\label{6.10}
\end{equation}
Next, we estimate $J_2$. We can write it in the form
$$
J_2=\intop_{\Omega^t}\bigg(\psi_{,rz}+{\psi_{,z}\over r}\bigg){1\over r}u_{,r}u_{,r}drdzdt'.
$$
Integrating by parts with respect to $r$ yields
$$\eqal{
J_2&=\intop_0^t\intop_{-a}^a\bigg(\psi_{,rz}+{\psi_{,z}\over r}\bigg){1\over r}u_{,r}u\bigg|_{r=0}^{r=R}dzdt'-\intop_{\Omega^t}\bigg(\psi_{,rz}+{\psi_{,z}\over r}\bigg)_{,r}{1\over r}u_{,r}udrdzdt'\cr
&\quad-\intop_{\Omega^t}\bigg(\psi_{,rz}+{\psi_{,z}\over r}\bigg)\,\bigg({1\over r}u_{,r}\bigg)_{,r}{u\over r}dxdt'\equiv J_{20}+J_{21}+J_{22},\cr}
$$
where the boundary term vanishes because $u|_{r=R}=0$ and (\ref{1.21}) implies
$$
\bigg(\psi_{,rz}+{\psi_{,z}\over r}\bigg){1\over r}u_{,r}u\bigg|_{r=0}=2d_{1,z}\bigg(v_{\varphi,r}+{v_\varphi\over r}\bigg)rv_\varphi\bigg|_{r=0}=4d_{1,z}b_1r^2b_1|_{r=0}=0.
$$
Using the transformation $\psi=r\psi_1$ in $J_{21}$ yields
$$\eqal{
J_{21}&=-\intop_{\Omega^t}(2\psi_{1,z}+r\psi_{1,rz})_{,r}{1\over r}{1\over r}u_{,r}udxdt'\cr
&=-\intop_{\Omega^t}(3\psi_{1,rz}+r\psi_{1,rrz}){1\over r}{1\over r}u_{,r}udxdt'.\cr}
$$
By the H\"older inequality, we have
$$\eqal{
|J_{21}|&\le 3\bigg|{\psi_{1,rz}\over r}\bigg|_{2,\Omega^t}\bigg|{1\over r}u_{,r}\bigg|_{2,\Omega^t}|u|_{\infty,\Omega^t}+|\psi_{1,rrz}|_{2,\Omega^t}\bigg|{1\over r}u_{1,r}\bigg|_{2,\Omega^t}|u|_{\infty,\Omega^t}\cr
&\le cD_1D_2|\Gamma_{,z}|_{2,\Omega^t},\cr}
$$
where (\ref{2.5}), (\ref{2.10}), (\ref{3.14}) and (\ref{3.15}) were used. Next, we consider $J_{22}$. Passing to function $\psi_1$, we get
$$
J_{22}=-\intop_{\Omega^t}(2\psi_{1,z}+r\psi_{1,rz}){1\over r}\bigg({1\over r}u_{,r}\bigg)_{,r}udxdt'.
$$
Hence
$$\eqal{
|J_{22}|&\le 2\bigg|\intop_{\Omega^t}{\psi_{1,z}\over r}\bigg({1\over r}u_{,r}\bigg)_{,r}udxdt'\bigg|+\bigg|\intop_{\Omega^t}\psi_{1,rz}\bigg({1\over r}u_{,r}\bigg)_{,r}udxdt'\bigg|\cr
&\equiv K_1+K_2,\cr}
$$
where $K_2$ is bounded by
$$
K_2\le\bigg|\intop_{\Omega^t}{\psi_{1,rz}\over r}u_{,rr}udxdt'\bigg|+\bigg|\intop_{\Omega^t}{\psi_{1,rz}\over r}{u_{,r}\over r}udxdt'\bigg|\equiv K_2^1+K_2^2.
$$
Using (\ref{2.5}), (\ref{2.10}), (\ref{5.2}) and (\ref{3.15}), we obtain
$$
|K_2^1|\le\bigg|{\psi_{1,rz}\over r}\bigg|_{2,\Omega^t}|u_{,rr}|_{2,\Omega^t}|u|_{\infty,\Omega^t}\le cD_2D_5|\Gamma_{,z}|_{2,\Omega^t}
$$
and
$$
|K_2^2|\le\bigg|{\psi_{1,rz}\over r}\bigg|_{2,\Omega^t}\bigg|{1\over r}u_{,r}\bigg|_{2,\Omega^t}|u|_{\infty,\Omega^t}\le cD_1D_2|\Gamma_{,z}|_{2,\Omega^t}.
$$
Summarizing,
\begin{equation}\eqal{
|J_2|&\le c|D_1D_2+D_2D_5]|\Gamma_{,z}|_{2,\Omega^t}+2\bigg|\intop_{\Omega^t}{\psi_{1,z}\over r}\bigg({1\over r}u_{,r}\bigg)_{,r}udxdt'\bigg|.\cr}
\label{6.11}
\end{equation}
Finally, we examine $J_3$. Using that $\psi=r\psi_1$ yields
$$
J_3=-\intop_{\Omega^t}(-r\psi_{1,zz}+3\psi_{1,r}+r\psi_{1,rr}){1\over r}u_{,r}{1\over r}u_{,z}dxdt'.
$$
Integrating by parts with respect to $z$, using that $\psi_1|_{S_2}=0$ and $\psi_{1,zz}|_{S_2}=-\Gamma|_{S_2}=0$, we obtain
$$\eqal{
J_3&=\intop_{\Omega^t}\bigg(-\psi_{1,zzz}+{3\over r}\psi_{1,rz}+\psi_{1,rrz}\bigg){1\over r}u_{,r}udxdt'\cr
&\quad+\intop_{\Omega^t}\bigg(-\psi_{1,zz}+{3\over r}\psi_{1,r}+\psi_{1,rr}\bigg)\,\bigg({1\over r}u_{,r}\bigg)_{,z}udxdt'\cr
&\equiv J_{31}+J_{32}.\cr}
$$
Using (\ref{2.5}), (\ref{2.10}), (\ref{3.14}) and (\ref{3.15}), we have
$$
|J_{31}|\le cD_1D_2|\Gamma_{,z}|_{2,\Omega^t}.
$$
To estimate $J_{32}$ we recall that
$$
-\psi_{1,rr}-{3\over r}\psi_{1,r}-\psi_{1,zz}=\Gamma.
$$
Then $J_{32}$ takes the form
$$
J_{32}=-\intop_{\Omega^t}(2\psi_{1,zz}+\Gamma)\bigg({1\over r}u_{,r}\bigg)_{,z}udxdt'\equiv J_{32}^1+J_{32}^2.
$$
Continuing,
$$\eqal{
|J_{32}^1|&\le c\bigg|\intop_{\Omega^t}{\psi_{1,zz}\over r^{1-\varepsilon_0}}u_{,rz}u^{1-\varepsilon_0}v_\varphi^{\varepsilon_0}dxdt\bigg|\cr
&\le cD_2^{1-\varepsilon_0}|v_\varphi|_{\infty,\Omega}^{\varepsilon_0}\bigg|{\psi_{1,zz}\over r^{1-\varepsilon_0}}\bigg|_{2,\Omega^t}\cdot|u_{,rz}|_{2,\Omega^t},\cr}
$$
where we used (\ref{2.10}). Finally, we have
$$
|J_{32}^1|\le cD_2^{1-\varepsilon_0}D_4{R^{\varepsilon_0}\over\varepsilon_0}|v_\varphi|_{\infty,\Omega^t}^{\varepsilon_0} |\Gamma_{,z}|_{2,\Omega^t}.
$$
Next,
$$\eqal{
|J_{32}^2|&\le\intop_{\Omega^t}\bigg|{\Gamma\over r^{1-\varepsilon_0}}\bigg||u_{,rz}|\,|u|^{1-\varepsilon_0}|v_\varphi|^{\varepsilon_0}dxdt'\cr
&\le cD_2^{1-\varepsilon_0}D_4{R^{\varepsilon_0}\over\varepsilon_0}|v_\varphi|_{\infty,\Omega^t}^{\varepsilon_0} |\Gamma_{,r}|_{2,\Omega^t}.\cr}
$$
Summarizing,
\begin{equation}
|J_3|\le cD_1D_2|\Gamma_{,z}|_{2,\Omega^t}+cD_2^{1-\varepsilon_0}D_4{R^{\varepsilon_0}\over\varepsilon_0}|v_\varphi|_{\infty,\Omega^t}^{\varepsilon_0} |\nabla\Gamma|_{2,\Omega^t}.
\label{6.12}
\end{equation}
Using estimates (\ref{6.10}), (\ref{6.11}), (\ref{6.12}) in (\ref{6.7}) implies
\begin{equation}\eqal{
|J|&\le\bigg|\intop_{\Omega^t}{\psi_{1,z}\over r^2}u_{,zz}udxdt\bigg|+2\bigg|\intop_{\Omega^t}{\psi_{1,z}\over r}\bigg({1\over r}u_{,r}\bigg)_{,r}udxdt\bigg|\cr
&\quad+c\bigg[{R^{\varepsilon_0}\over\varepsilon_0}(D_1+D_4)D_2^{1-\varepsilon_0} |v_\varphi|_{\infty,\Omega^t}^{\varepsilon_0}\cr
&\quad+D_1D_2+D_2D_4+D_2D_5\bigg]\,|\nabla\Gamma|_{2,\Omega^t}.\cr}
\label{6.13}
\end{equation}
Using estimates (\ref{6.4}), (\ref{6.5}), (\ref{6.13}) in (\ref{6.3}) implies the inequality
\begin{equation}\eqal{
&|\omega_r(t)|_{2,\Omega}^2+|\omega_z(t)|_{2,\Omega}^2+\nu(|\nabla\omega_r|_{2,\Omega^t}^2+ |\nabla\omega_z|_{2,\Omega^t}^2+|\Phi|_{2,\Omega^t}^2)\cr
&\le\bigg|\intop_{\Omega^t}{\psi_{1,z}\over r^2}u_{,zz}udxdt'\bigg|+2\bigg|\intop_{\Omega^t}{\psi_{1,z}\over r}\bigg({1\over r}u_{,r}\bigg)_{,r}udxdt'\bigg|\cr
&\quad+c\bigg[{R^{\varepsilon_0}\over\varepsilon_0}(D_1+D_4)D_2^{1-\varepsilon_0} |v_\varphi|_{\infty,\Omega^t}^{\varepsilon_0}+D_1D_2+D_2(D_4+D_5)\bigg]\cdot\cr
&\quad\cdot|\nabla\Gamma|_{2,\Omega^t}+c|f_\varphi|_{2,S_1^t}(D_4+D_5)\cr
&\quad+c(|F_r|_{6/5,2,\Omega^t}^2+|F_z|_{6/5,2,\Omega^t}^2)+|\omega_r(0)|_{2,\Omega}^2+ |\omega_z(0)|_{2,\Omega}^2.\cr}
\label{6.14}
\end{equation}
Recalling that $\omega_r=-{1\over r}u_{,z}$, $\omega_z={1\over r}u_{,r}$ (see (\ref{1.5})) and applying the H\"older and Young inequalities to the first term on the r.h.s. of (\ref{6.14}) we bound it by
$$
\varepsilon_1|\omega_{r,z}|_{2,\Omega^t}^2+{1\over 4\varepsilon_1}\intop_{\Omega^t}{\psi_{1,z}^2\over r^2}u^2dxdt'\equiv L_1.
$$
The second term in $L_1$ can be written in the form
$$
\intop_{\Omega^t}{\psi_{1,z}^2\over r^{2(1-\varepsilon_0)}}{u^2\over r^{2\varepsilon_0}}dxdt'\le D_2^{2(1-\varepsilon_0)}|v_\varphi|_{\infty,\Omega^t}^{2\varepsilon_0}\intop_{\Omega^t} {\psi_{1,z}^2\over r^{2(1-\varepsilon_0)}}dxdt',
$$
where $\varepsilon_0$ can be chosen as small as we need. By the Hardy inequality
$$
\intop_{\Omega^t}{\psi_{1,z}^2\over r^{2(1-\varepsilon_0)}}dxdt'\le{R^{2\varepsilon_0}\over\varepsilon_0^2}\intop_{\Omega^t}\psi_{1,rz}^2dxdt'.
$$
Applying the interpolation inequality (\ref{2.3}) (see [Ch. 3, sect. 15]\cite{BIN})
$$
\intop_\Omega\psi_{1,rz}^2dx\le\bigg(\intop_\Omega|\nabla^2\psi_{1,z}|^2dx\bigg)^\theta \bigg(\intop_\Omega\psi_{1,z}^2dx\bigg)^{1-\theta},
$$
where $\theta$ satisfies the equality
$$
{3\over 2}-1=(1-\theta){3\over 2}+\theta\bigg({3\over 2}-2\bigg)\quad {\rm so}\ \ \theta=1/2.
$$
Using (\ref{2.12}), we get
$$
\intop_{\Omega^t}\psi_{1,zr}^2dx\le c|\nabla^2\psi_{1,z}|_{2,\Omega}|\psi_{1,z}|_{2,\Omega}\le cD_1|\nabla^2\psi_{1,z}|_{2,\Omega}.
$$
Summarizing,
\begin{equation}\eqal{
&\bigg|\intop_{\Omega^t}{\psi_{1,z}\over r^2}u_{,zz}udxdt\bigg|\le\varepsilon_1|\omega_{r,z}|_{2,\Omega^t}^2\cr
&\quad+{c\over 4\varepsilon_1}D_1D_2^{2(1-\varepsilon_0)}{R^{2\varepsilon_0}\over\varepsilon_0^2} |v_\varphi|_{\infty,\Omega^t}^{2\varepsilon_0}|\Gamma_{,z}|_{2,\Omega^t}.\cr}
\label{6.15}
\end{equation}
Similarly, the second term on the r.h.s. of (\ref{6.14}) is bounded by
\begin{equation}\eqal{
&\bigg|\intop_{\Omega^t}{\psi_{1,z}\over r}\bigg({1\over r}u_{,r}\bigg)_{,r}udxdt'\bigg|\le\varepsilon_2|\omega_{z,r}|_{2,\Omega^t}^2\cr
&\quad+{c\over 4\varepsilon_2}D_1D_2^{2(1-\varepsilon_0)}{R^{2\varepsilon_0}\over\varepsilon_0^2} |v_\varphi|_{\infty,\Omega^t}^{2\varepsilon_0}|\Gamma_{,z}|_{2,\Omega^t}.\cr}
\label{6.16}
\end{equation}
Using (\ref{6.15}) and (\ref{6.16}) in (\ref{6.14}) yields (\ref{6.1}).
\end{proof}

\begin{lemma}\label{l6.2}
Assume that $D_1$, $D_2$ are defined in Notation \ref{n1.1}, $v_\varphi(0)\in L_\infty(\Omega)$, $f_\varphi/r\in L_1(0,t;L_\infty(\Omega))$. Then
\begin{equation}
|v_\varphi(t)|_{\infty,\Omega}\le{D_2\over\sqrt{\nu}}D_1^{1/4}X^{3/4}+D_7,
\label{6.17}
\end{equation}
where
\begin{equation}
D_7=\sqrt{2}D_2^{1/2}\bigg|{f_\varphi\over r}\bigg|_{\infty,1,\Omega^t}^{1/2}+|v_\varphi(0)|_{\infty,\Omega}.
\label{6.18}
\end{equation}
\end{lemma}

\begin{proof}
Multiplying $(\ref{1.7})_2$ by $v_\varphi|v_\varphi|^{s-2}$ and integrating over $\Omega$ yields
\begin{equation}\eqal{
&{1\over s}{d\over dt}|v_\varphi|_{s,\Omega}^s+{4\nu(s-1)\over s^2}|\nabla|v_\varphi|^{s/2}|_{2,\Omega}^2+\nu\intop_\Omega{|v_\varphi|^s\over r^2}dx\cr
&=\intop_\Omega\psi_{1,z}|v_\varphi|^sdx+\intop_\Omega f_\varphi v_\varphi^{s-1}dx\cr}
\label{6.19}
\end{equation}
The first term on the r.h.s. of (\ref{6.19}) is estimated by
$$
\varepsilon\intop_\Omega{|v_\varphi|^s\over r^2}dx+{D_2^2\over 4\varepsilon}\intop_\Omega|\psi_{1,z}|^2|v_\varphi|^{s-2}dx
$$
The second integral on the r.h.s. of \eqref{6.10} is estimated by
$$\eqal{
&\intop_\Omega|f_\varphi|\,|v_\varphi|^{s-1}dx=\intop_\Omega\bigg|{f_\varphi\over r}\bigg|r|v_\varphi|^{s-1}dx\cr
&\le D_2\intop_\Omega\bigg|{f_\varphi\over r}\bigg|\,|v_\varphi|^{s-2}dx\le D_2\bigg|{f_\varphi\over r}\bigg|_{s/2,\Omega}|v_\varphi|_{s,\Omega}^{s-2}.\cr}
$$
In view of the above estimates inequality (\ref{6.19}) reads
$$
{1\over s}{d\over dt}|v_\varphi|_{s,\Omega}^s\le{D_2^2\over 2\nu}|\psi_{1,z}|_{s,\Omega}^2|v_\varphi|_{s,\Omega}^{s-2}+D_2\bigg|{f_\varphi\over r}\bigg|_{s/2,\Omega}|v_\varphi|_{s,\Omega}^{s-2}.
$$
Simplifying, we get
$$
{d\over dt}|v_\varphi|_{s,\Omega}^2\le{D_2^2\over\nu}|\psi_{1,z}|_{s,\Omega}^2+2D_2\bigg|{f_\varphi\over r}\bigg|_{s/2,\Omega}.
$$
Integrating the inequality with respect to time and passing with $s$ to $\infty$ we get
$$
|v_\varphi|_{\infty,\Omega}^2\le{D_2^2\over\nu}\intop_0^t|\psi_{1,z}|_{\infty,\Omega}^2dt'+2D_2\bigg|{f_\varphi\over r}\bigg|_{\infty,1,\Omega^t}+|v_\varphi(0)|_{\infty,\Omega}^2.
$$
Using the interpolation
$$
|\psi_{1,z}|_{\infty,\Omega}\le|\psi_{1,z}|_{2,\Omega}^{1/4}|D^2\psi_{1,z}|_{2,\Omega}^{3/4}
$$
and (\ref{2.3}), (\ref{3.14}), we obtain
$$
|v_\varphi|_{\infty,\Omega}^2\le{D_2^2\over\nu}D_1^{1/2}|\Gamma_{,z}|_{2,\Omega^t}^{3/2}+\bigg|{f_\varphi\over r}\bigg|_{\infty,1,\Omega^t}+2D_2|v_\varphi(0)|_{\infty,\Omega}^2.
$$
The above inequality implies (\ref{6.17}) and concludes the proof.
\end{proof}

\section{Global estimate for regular solutions}\label{s7}

Assuming appropriate regularity of data we show that boundedness of $X(t)$ implies that
\begin{equation}
\|v\|_{W_2^{4,2}(\Omega^t)}+\|\nabla p\|_{W_2^{2,1}(\Omega^t)}\le C,
\label{7.1}
\end{equation}
where $C>0$ depends only on initial data and the forcing.

\subsection{Preliminaries}\label{s7.1}

We first introduce some functional analytic tools to handle the anisotropic Sobolev spaces.

\begin{definition}[Anisotropic Sobolev and Sobolev-Slobodetskii spaces]\label{d7.1}
We denote by
\begin{itemize}
\item[1.] $W_{p,p_0}^{k,k/2}(\Omega^t)$, $k,k/2\in\N\cup\{0\}$, $p,p_0\in[1,\infty]$ -- the anisotropic Sobolev space with a mixed norm, which is a completion of $C^\infty(\Omega^T)$-functions under the norm
$$
\|u\|_{W_{p,p_0}^{k,k/2}(\Omega^T)}=\bigg(\intop_0^T\bigg(\sum_{|\alpha|+2a\le k}\intop_\Omega|D_x^\alpha\partial_t^au|^p\bigg)^{p_0/p}dt\bigg)^{1/p_0}.
$$
\item[2.] $W_{p,p_0}^{s,s/2}(\Omega^T)$, $s\in\R_+$, $p,p_0\in[1,\infty)$ -- the Sobolev-Slobodetskii space with the finite norm 
$$\eqal{
&\|u\|_{W_{p,p_0}^{s,s/2}(\Omega^T)}=\sum_{|\alpha|+2a\le|s|}\|D_x^\alpha\partial_t^au\|_{L_{p,p_0}(\Omega^T)}\cr
&+\bigg[\intop_0^T\!\!\bigg(\intop_\Omega\intop_\Omega\!\sum_{|\alpha|+2a\le[s]}\!\! {|D_x^\alpha\partial_t^au(x,t)-D_{x'}^\alpha\partial_t^au(x',t)|^p\over|x-x'|^{n+p(s-[s])}}dxdx'\bigg)^{p_0/p}dt \bigg]^{1/p_0}\cr
&+\bigg[\intop_\Omega\!\!\bigg(\intop_0^T\intop_0^T\!\sum_{|\alpha|+2a=[s]}\!\! {|D_x^\alpha\partial_t^au(x,t)-D_x^\alpha\partial_{t'}^au(x,t')|^{p_0}\over|t-t'|^{1+p_0({s\over 2}-[{s\over 2}])}}dtdt'\bigg)^{p/p_0}dx\bigg]^{1/p}\!,\cr}
$$
where $a\in\N\cup\{0\}$, $[s]$ is the integer part of $s$ and $D_x^\alpha$ denotes the partial derivative in the spatial variable $x$ cooresponding to multiindex $\alpha$. For $s$ odd the last but one term in the above norm vanishes whereas for $s$ even the last two terms vanish. We also use notation $L_p(\Omega^T)=L_{p,p}(\Omega^T)$, $W_p^{s,s/2}(\Omega^T)=W_{p,p}^{s,s/2}(\Omega^T)$.
\item[3.] $B_{p,p_0}^l(\Omega)$, $l\in\R_+$, $p,p_0\in[1,\infty)$ -- the Besov space with the finite norm
$$
\|u\|_{B_{p,p_0}^l(\Omega)}=\|u\|_{L_p(\Omega)}+\bigg(\sum_{i=1}^n\intop_0^\infty {\|\Delta_i^m(h,\Omega)\partial_{x_i}^ku\|_{L_p(\Omega)}^{p_0}\over h^{1+(l-k)_{p_0}}}dh\bigg)^{1/p_0},
$$
where $k\in\N\cup\{0\}$, $m\in\N$, $m>l-k>0$, $\Delta_i^j(h,\Omega)u$, $j\in\N$, $h\in\R_+$ is the finite difference of the order $j$ of the function $u(x)$ with respect to $x_i$ with
$$\eqal{
&\Delta_i^1(h,\Omega)u=\Delta_i(h,\Omega)\cr
&=u(x_1,\dots,x_{i-1},x_i+h,x_{i+1},\dots,x_n)-u(x_1,\dots,x_n),\cr
&\Delta_i^j(h,\Omega)=\Delta_i(h,\Omega)\Delta_i^{j-1}(h,\Omega)u\quad {\rm and}\ \ \Delta_i^j(h,\Omega)u=0\cr
&{\rm for}\ \ x+jh\not\in\Omega.\cr}
$$
In has been proved in \cite{G} that the norms of the Besov space $B_{p,p_0}^l(\Omega)$ are equivalent for different $m$ and $k$ satisfying the condition $m>l-k>0$.
\end{itemize}
\end{definition}

We need the following interpolation lemma.

\begin{lemma}[Anisotropic interpolation, see {\cite[Ch. 4, Sect. 18]{BIN}}]\label{l7.2}
Let $u\in W_{p,p_0}^{s,s/2}(\Omega^T)$, $s\in\R_+$, $p,p_0\in[1,\infty]$, $\Omega\subset\R^3$. Let $\sigma\in\R_+\cup\{0\}$, and
$$
\varkappa={3\over p}+{2\over p_0}-{3\over q}-{2\over q_0}+|\alpha|+2a+\sigma<s.
$$
Then $D_x^\alpha\partial_t^au\in W_{q,q_0}^{\sigma,\sigma/2}(\Omega^T)$, $q\ge p$, $q_0\ge p_0$ and there exists $\varepsilon\in(0,1)$ such that
$$
\|D_x^\alpha\partial_t^au\|_{W_{q,q_0}^{\sigma,\sigma/2}(\Omega^T)}\le\varepsilon^{s-\varkappa} \|u\|_{W_{p,p_0}^{s,s/2}(\Omega^t)}+c\varepsilon^{-\varkappa}\|u\|_{L_{p,p_0}(\Omega^t)}.
$$
We recall from \cite{B} the trace and the inverse trace theorems for Sobolev spaces with a mixed norm.
\end{lemma}

\begin{lemma}\label{l7.3}
(traces in $W_{p,p_0}^{s,s/2}(\Omega^T)$, see \cite{B})
\begin{itemize}
\item[(i)] Let $u\in W_{p,p_0}^{s,s/2}(\Omega^t)$, $s\in\R_+$, $p,p_0\in(1,\infty)$. Then $u(x,t_0)=u(x,t)|_{t=t_0}$ for $t_0\in[0,T]$ belongs to $B_{p,p_0}^{s-2/p_0}(\Omega)$, and
$$
\|u(\cdot,t_0)\|_{B_{p,p_0}^{s-2/p_0}(\Omega)}\le c\|u\|_{W_{p,p_0}^{s,s/2}(\Omega^T)},
$$
where $c$ does not depend on $u$.
\item[(ii)] For given $\bar u\in B_{p,p_0}^{s-2/p_0}(\Omega)$, $s\in\R_+$, $s>2/p_0$, $p_0\in(1,\infty)$, there exists a function $u\in W_{p,p_0}^{s,s/2}(\Omega^t)$ such that $u|_{t=t_0}=\bar u$ for $t_0\in[0,T]$ and
$$
\|u\|_{W_{p,p_0}^{s,s/2}(\Omega^T)}\le c\|\bar u\|_{B_{p,p_0}^{s-2/p_0}(\Omega)},
$$
where constant $c$ does not depend on $\bar u$.
\end{itemize}
\end{lemma}

We need the following imbeddings between Besov spaces

\begin{lemma}[see {\cite[Th. 4.6.1]{T}}]\label{l7.4}
Let $\Omega\subset\R^n$ be an arbitrary domain.
\begin{itemize}
\item[(a)] Let $s\in\R_+$, $\varepsilon>0$, $p\in(1,\infty)$, and $1\le q_1\le q_2\le\infty$. Then
$$
B_{p,1}^{s+\varepsilon}(\Omega)\subset B_{p,\infty}^{s+\varepsilon}(\Omega)\subset B_{p,q_1}^s(\Omega)\subset B_{p,q_2}^s(\Omega)\subset B_{p,1}^{s-\varepsilon}(\Omega)\subset B_{p,\infty}^{s-\varepsilon}(\Omega).
$$
\item[(b)] Let $\infty>q\ge p>1$, $1\le r\le\infty$, $0\le t\le s<\infty$ and
$$
t+{n\over p}-{n\over q}\le s.
$$
then $B_{p,r}^s(\Omega)\subset B_{q,r}^t(\Omega)$.
\end{itemize}
\end{lemma}

\begin{lemma}[see {\cite[Ch. 4, Th. 18.8]{BIN}}]\label{lemma 7.5}
Let $1\le\theta_1<\theta_2\le\infty$. Then
$$
\|u\|_{B_{p,\theta_2}^l(\Omega)}\le c\|u\|_{B_{p,\theta_1}^l(\Omega)},
$$
where $c$ does not depend on $u$.
\end{lemma}

\begin{lemma}[see {\cite[Ch. 4, Th. 18.9]{BIN}}]\label{l7.6}
Let $l\in\N$ and $\Omega$ satisfy the $l$-horn condition. Then the following imbeddings hold
$$\eqal{
&\|u\|_{B_{p,2}^l(\Omega)}\le c\|u\|_{W_p^l(\Omega)}\le c\|u\|_{B_{p,p}^l(\Omega)},\quad &1\le p\le 2,\cr
&\|u\|_{B_{p,p}^l(\Omega)}\le c\|u\|_{W_p^l(\Omega)}\le c\|u\|_{B_{p,2}^l(\Omega)},\quad &2\le p<\infty,\cr
&\|u\|_{B_{p,\infty}^l(\Omega)}\le c\|u\|_{W_p^l(\Omega)}\le c\|u\|_{B_{p,1}^l(\Omega)},\quad &1\le p\le\infty.\cr}
$$
\end{lemma}

Consider the nonstationary Stokes system in $\Omega\subset\R^3$:
$$
\begin{aligned}
    &v_t-\nu\Delta v+\nabla p=f,\\
&\divv v=0, \\
&v_r=v_\varphi=\omega_\varphi=0&&\text{ on }S_1^T,\\
&v_z=\omega_\varphi=v_{\varphi,z}&&\text{ on }S_2^T,\\
&v|_{t=0}=v(0),
\end{aligned}
$$
with the boundary conditions (\ref{1.2}) and given initial condition $v(0)$.

\begin{lemma}[see \cite{MS}]\label{l7.7}
Assume that $f\in L_{q,r}(\Omega^T)$, $v(0)\in B_{q,r}^{2-2/r}(\Omega)$, $r,q\in(1,\infty)$. Then there exists a unique solution to the above system such that $v\in W_{q,r}^{2,1}(\Omega^T)$, $\nabla p\in L_{q,r}(\Omega^T)$ with the following estimate
\begin{equation}\eqal{
\|v\|_{W_{q,r}^{2,1}(\Omega^T)}+\|\nabla p\|_{L_{q,r}(\Omega^t)}&\le  c(\|f\|_{L_{q,r}(\Omega^t)}+\|v(0)\|_{B_{r,q}^{2-2/r}(\Omega)}).\cr}
\label{7.2}
\end{equation}
\end{lemma}

\subsection{Proof of (\ref{7.1})}\label{s7.2}
We show (\ref{7.1}) in the following series of lemmas.

\begin{lemma}\label{l7.8}
Suppose that $X(t)<\infty$, i.e. that
\begin{equation}
\|\Phi\|_{V(\Omega^t)}+\|\Gamma\|_{V(\Omega^t)}\le\phi_1,
\label{7.3}
\end{equation}
where $\phi_1$ depends only on the initial data and the forcing. Assume that
$$
f\in W_2^{2,1}(\Omega^t),\quad v(0)\in W_2^3(\Omega).
$$
Then
\begin{equation}
\|v\|_{W_2^{4,2}(\Omega^t)}+\|\nabla p\|_{W_2^{2,1}(\Omega^t)}\le\phi(\phi_1,\|f\|_{W_2^{2,1}(\Omega^t)},\|v(0)\|_{H^3(\Omega)}).
\label{7.4}
\end{equation}
\end{lemma}

\begin{proof}
From (\ref{7.3}) we have
\begin{equation}
\|\Gamma\|_{V(\Omega^t)}\le\phi_1.
\label{7.5}
\end{equation}
Section 3.2 implies
\begin{equation}
\|\psi_1\|_{2,\infty,\Omega^t}\le c\phi_1,
\label{7.6}
\end{equation}
\begin{equation}
\|\psi_1\|_{3,2,\Omega^t}\le c\phi_1.
\label{7.7}
\end{equation}
From (\ref{1.19}) the following relations hold
\begin{equation}
v_r=-r\psi_{1,z},\quad v_z=2\psi_1+r\psi_{1,r}.
\label{7.8}
\end{equation}
Hence (\ref{7.6}) and $R$ finite imply
\begin{equation}
\|v_r\|_{1,2,\infty,\Omega^t}+\|v_z\|_{1,2,\infty,\Omega^t}\le c\phi_1.
\label{7.9}
\end{equation}
To increase regularity of $v$ we consider the Stokes problem
\begin{equation}\eqal{
&v_{,t}-\nu\Delta v+\nabla p=-v'\cdot\nabla v+f\quad &{\rm in}\ \ \Omega^T,\cr
&\divv v=0\quad &{\rm in}\ \ \Omega^T,\cr
&v\cdot\bar n=0,\ \ v_{z,r}=0,\ \ v_\varphi=0\quad &{\rm on}\ \ S_1^T,\cr
&v\cdot\bar n=0,\ \ v_{r,z}=0,\ \ v_{\varphi,z}=0\quad &{\rm on}\ \ S_2^T,\cr
&v|_{t=0}=v(0)\quad &{\rm in}\ \ \Omega,\cr}
\label{7.10}
\end{equation}
where $v'=v_r\bar e_r+v_z\bar e_z$.

From (\ref{7.9}) we have
\begin{equation}
|v'|_{6,\infty,\Omega^t}\le c\phi_1.
\label{7.11}
\end{equation}
For solutions to (\ref{7.10}) the following energy estimate holds
\begin{equation}
\|v\|_{V(\Omega^t)}\le c(|f|_{2,\Omega^t}+|v(0)|_{2,\Omega})\equiv d_1.
\label{7.12}
\end{equation}
Hence, (\ref{7.12}) yields
\begin{equation}
|\nabla v|_{2,\Omega^t}\le d_1.
\label{7.13}
\end{equation}
Estimates (\ref{7.11}) and (\ref{7.13}) imply
\begin{equation}
|v'\cdot\nabla v|_{3/2,2,\Omega^t}\le c\phi_1d_1.
\label{7.14}
\end{equation}
Applying \cite{MS} to (\ref{7.10}) yields
\begin{equation}\eqal{
&\|v\|_{W_{{3\over 2},2}^{2,1}(\Omega^t)}+|\nabla p|_{3/2,2,\Omega^t}\cr
&\le c(|f|_{3/2,2,\Omega^t}+\|v(0)\|_{B_{3/2,2}^1(\Omega)}+\psi_1d_1)\equiv d_2.\cr}
\label{7.15}
\end{equation}
In view of the imbedding (see \cite[Ch. 3, Sect. 10]{BIN})
\begin{equation}
|\nabla v|_{5/2,\Omega^t}\le c\|v\|_{W_{3/2,2}^{2,1}(\Omega^t)}
\label{7.16}
\end{equation}
and (\ref{7.11}) we derive that
\begin{equation}
|v'\cdot\nabla v|_{{30\over 17},{5\over 2},\Omega^t}\le c\phi_1d_2.
\label{7.17}
\end{equation}
Then applying again \cite{MS} to problem (\ref{7.10}) yields
\begin{equation}\eqal{
&\|v\|_{W_{{30\over 17},{5\over 2}}^{2,1}(\Omega^t)}+|\nabla p|_{{30\over 17},{5\over 2},\Omega^t}\cr
&\le c(|f|_{{30\over 17},{5\over 2},\Omega^t}+\|v(0)\|_{B_{{30\over 17},{5\over 2}}^{2-4/5}(\Omega)}+\phi_1d_2)\equiv d_3.\cr}
\label{7.18}
\end{equation}
In view of the imbedding (see \cite[Ch. 3, Sect. 10]{BIN})
\begin{equation}
|\nabla v|_{{10\over 3},\Omega^t}\le c\|v\|_{W_{{30\over 17},{5\over 2}}^{2,1}(\Omega^t)}
\label{7.19}
\end{equation}
and (\ref{7.11}) we have
\begin{equation}
|v'\cdot\nabla v|_{{15\over 7},{10\over 3},\Omega^t}\le c\phi_1d_3.
\label{7.20}
\end{equation}
Applying \cite{MS} to (\ref{7.10}) implies
\begin{equation}\eqal{
&\|v\|_{W_{{15\over 7},{10\over 3}}^{2,1}(\Omega^t)}+|\nabla p|_{{15\over 7},{10\over 3},\Omega^t}\cr
&\le c(|f|_{{15\over 7},{10\over 3},\Omega^t}+\|v(0)\|_{B_{{15\over 7},{10\over 3}}^{2-6/10}(\Omega)}+\phi_1d_3)=d_4.\cr}
\label{7.21}
\end{equation}
Lemma \ref{l7.3} yields
\begin{equation}
\|v\|_{L_\infty(0,t;B_{{15\over 7},{10\over 3}}^{2-6/10}(\Omega^t))}\le c\|v\|_{W_{{15\over 7},{10\over 3}}^{2,1}(\Omega^t)}.
\label{7.22}
\end{equation}
Theorem 18.10 from \cite{BIN} gives
\begin{equation}
|v(t)|_{q,\Omega}\le c\|v\|_{B_{{15\over 7},{10\over 3}}^{7/5}(\Omega)}.
\label{7.23}
\end{equation}
The estimate holds for any finite $q$ because it satisfies the relation $7/5\ge 7/5-3/q$.

Next, we use the imbedding (see \cite[Ch. 3, Sect. 10]{BIN})
\begin{equation}
|\nabla v|_{5,\Omega^t}\le c\|v\|_{W_{{15\over 7},{10\over 3}}^{2,1}(\Omega^t)}.
\label{7.24}
\end{equation}
From (\ref{7.23}) and (\ref{7.24}) we have
\begin{equation}
|v\cdot\nabla v|_{5',\Omega^t}\le cd_4^2,
\label{7.25}
\end{equation}
where $5'<5$ but it is arbitrary close to 5.

In view of (\ref{7.25}) and \cite{MS} we have
\begin{equation}
\|v\|_{W_{5'}^{2,1}(\Omega^t)}+|\nabla p|_{5',\Omega^t}\le c(|f|_{5',\Omega^t}+\|v(0)\|_{W_{5'}^{2-2/5'}(\Omega)}+d_4^2)\equiv d_5.
\label{7.26}
\end{equation}
From (\ref{7.26}) it follows that $v\in L_\infty(\Omega^t)$ and $\nabla v\in L_q(\Omega^t)$ for any finite $q$.

Then
$$
|\nabla(v'\cdot\nabla v)|_{5',\Omega^t}\le cd_5^2
$$
and $$
|\partial_t^{1/2}(v'\cdot\nabla v)|_{5',\Omega^t}\le cd_5^2,
$$
where $\partial_t^{1/2}$ denotes the fractional partial derivative in time.

Then \cite{MS} implies
\begin{equation}\eqal{
&\|v\|_{W_{10/3}^{3,3/2}(\Omega^t)}+\|\nabla p\|_{W_{10/3}^{1,1/2}(\Omega^t)}\cr
&\le c(\|f\|_{W_{10/3}^{1,1/2}(\Omega^t)}+\|v(0)\|_{W_{10/3}^{3-2/5'}(\Omega)}+d_5^2)\equiv d_6.\cr}
\label{7.27}
\end{equation}
Continuing the considerations yields
\begin{equation}
\|v\|_{W_2^{4,2}(\Omega^t)}+\|\nabla p\|_{W_2^{2,1}(\Omega^t)}\le c(\|f\|_{W_2^{2,1}(\Omega^t)}+\|v(0)\|_{W_2^3(\Omega)}+d_6^2).
\label{7.28}
\end{equation}
This implies (\ref{7.4}) and ends the proof.
\end{proof}

\noindent
{\bf Conflict of interest statement}

The authors report there are no competing interests to declare.

\noindent
{\bf Data availability statement}

The authors report that there is no data associated with this work.

\bibliographystyle{amsplain}

\end{document}